\documentclass[reqno]{amsart}
\usepackage{graphicx} 
\usepackage{subfiles}
\usepackage{enumitem}
\usepackage[T1]{fontenc}
\usepackage{euscript}
\usepackage{dsfont}
\newtheorem{theorem}{Theorem}[section]
\newtheorem{lemma}{Lemma}[section]

\theoremstyle{remark}
\newtheorem{remark}{Remark}[section]
\newcommand*\diff{\mathop{}\!\mathrm{d}}
\usepackage{physics}
\newcommand{\R}{\mathbb{R}}
\newcommand{\eps}{\epsilon}
\newcommand{\N}{\mathbb{N}}
\title{Viscous Destabilization for Large Shocks of Conservation Laws}
\author{Paul Blochas}
\address{Department of Mathematics, The University of Texas at Austin, Austin, TX 78712.}
\email{paul.blochas@utexas.edu}

\author{Jeffrey Cheng}
\address{Department of Mathematics, The University of Texas at Austin, Austin, TX 78712.}
\email{jeffrey.cheng@utexas.edu}

\date{\today}
\thanks{\textbf{Acknowledgement}: The authors would like to thank their mentor Alexis Vasseur for the suggestion of this problem and guidance throughout the project. J. Cheng is partially supported by NSF grant: DMS-1840314.}

\begin{document}
\begin{abstract}
The recent theory of $a-$contraction with shifts provides \newline $L^2$-stability for shock waves of $1-$D hyperbolic systems of conservation laws. The theory has been established at the inviscid level uniformly in the shock amplitude, and at the viscous level for small shocks. In this work, we investigate whether the $a-$contraction property holds uniformly in the shock amplitude for some specific systems with viscosity. We show that in some cases, the $a-$contraction fails for sufficiently large shocks. This showcases a ``viscous destabilization'' effect in the sense that the $a$-contraction property is verified for the inviscid model, but can fail for the viscous one. This also shows that the $a$-contraction property, even among small perturbations, is stronger than the classical notion of nonlinear stability, which is known to hold regardless of shock amplitude for viscous scalar conservation laws. 
\end{abstract}
\maketitle
\tableofcontents

\section{Introduction \& main results}
\subsection{Introduction}
\par In this paper, we will discuss stability properties of some $1-$D systems of viscous conservation laws. \\
The general form of such systems is as follows
\begin{align*}
    u_t+(A(u))_x=(B(u)u_x)_x.
\end{align*}
Here, we will consider two cases, a scalar case and a special barotropic Navier-Stokes system of equations. \\
First, a semilinear scalar case:
\begin{equation}\label{viscous}
\begin{cases}
u_t+(A(u))_x=u_{xx}, & (x,t) \in \mathbb{R} \times (0,\infty), \\ 
u(x,0)=u_0(x), & x \in \mathbb{R}.
\end{cases}
\end{equation}
Given any smooth, strictly convex flux $A$ and any $u_- > u_+$, there exists a smooth function $S$, and $\sigma \in \mathbb{R}$, such that $S(x-\sigma t)$ is a ``(viscous) shock'' solution to
\begin{equation}\label{viscousshock}
-\sigma S'(\xi)+A(S)'(\xi)=S''(\xi), 
\end{equation}
connecting two distinct end states $u_-$ and $u_+$, with speed $\sigma=\frac{A(u_+)-A(u_-)}{u_+-u_-}$. Integrating the equation gives the following: $-\sigma (S-u_{\pm})+A(S)-A(u_{\pm})=S'$.

\subsubsection{Previous work}
\par There have been many past studies on the stability of shocks of viscous conservation laws. Long time behavior for Burgers' equation was studied by Hopf \cite{MR0047234} and Nishihara \cite{MR0839318}. For more general fluxes, results for small perturbations of the viscous shock wave $S$ were given in the scalar case by Il'in and Oleinik \cite{MR0120469}. This problem was also considered using energy methods for zero mass perturbations in \cite{MR0853782}. Goodman generalized the zero mass perturbation method to some multidimensional problems \cite {MR0978372}. Sattinger developed a semigroup strategy which allowed him to obtain stability in weighted spaces \cite{MR0435602}. Szepessy and Xin on one hand \cite{MR1207241} and Liu \cite{MR0791863}, \cite{MR1470318} on the other hand expanded these results to the system case, both by the use of sophisticated energy methods. Howard in \cite{MR1699970} revisited the problem with pointwise Green's function methods. Zumbrun, with Howard and other coauthors, developed this method further, which, combined with high-frequency estimates, extended the results to multidimensional systems (see \cite{MR1362168} and references cited there). Freist\"uhler and Serre showed $L^1$-stability without the smallness condition by using techniques based on the maximum principle \cite{MR1488516} for scalar one-dimensional problems.

\subsubsection{The relative entropy method}
\par This paper regards recent developments of the relative entropy method initiated by Dafermos and DiPerna in 1979, who proved a weak-strong stability theorem for the inviscid model in \cite{MR0546634}, \cite{MR0523630}:
\begin{equation}\label{inviscid}
\begin{cases}
u_t+(A(u))_x=0, & (x,t) \in \mathbb{R} \times (0,\infty), \\
u(x,0)=u_0(x), & x \in \mathbb{R}. \\
\end{cases}
\end{equation}
\par This line of work shows stability independent of the size of the initial perturbation. Leger proved relative entropy stability up to shift for shocks for the inviscid model \cite{MR2771666}. That is, there exists a Lipschitz function $X:[0,\infty) \to \mathbb{R}$ such that for every $t$ positive:
\[
\int_\mathbb{R}\eta(u(x+X(t),t)|S(x+\sigma t))\diff x \leq \int_\mathbb{R}\eta(u(x,0)|S(x))\diff x,
\]
where $u$ is a solution to \eqref{inviscid} and $S$ the shock of reference, both entropic for $\eta$. Here $\eta(u|v)=\eta(u)-\eta(v)-\eta'(v)(u-v)$ is the relative entropy associated to the entropy $\eta$. In the scalar case, by picking $\eta(u)=u^2$, the above result shows exactly a contraction in $L^2$ up to a suitable shift. However, it has also been shown that for systems of conservation laws, the $L^2$-contraction need not hold (even up to shift) \cite{MR3322780}. This issue was circumvented for extremal shocks. It was proved that in such cases, the following modification of this notion of stability, called $a-$contraction, holds.

\subsubsection{The theory of $a$-contraction}
\par The theory of $a$-contraction up to shift focuses on showing that, by weighting the relative entropy adequately, the following quantity is non-increasing in $t$:
\begin{align*}
    \int_{\R} a(x)\eta(u(t,x+X(t))|S(x))\diff x.
\end{align*}
By picking $a$ uniformly bounded away from $0$ and $+\infty$, the $a-$contraction is seen to imply $L^2$-stability for solutions taking values in a compact subspace of the domain of the entropy (as the integral of the relative entropy itself is equivalent to the square of the $L^2$ distance in that case). In the inviscid setting, $a$-contraction has been proven for extremal shocks of systems of conservation laws \cite{MR3519973}, \cite{MR3537479}. 
\par In regards to the viscous theory, Kang and Vasseur showed that the $L^2$-contraction up to shift holds when $A$ is a small perturbation of the Burgers' flux \cite{MR3592682}. However, they also show that even in the scalar setting, for a fixed shock size, there are fluxes for which the $L^2$-contraction (i.e. $a$-contraction with $a=1$) fails to hold. Due to this, Stokols proved $a$-contraction for more general fluxes with a non-trivial weight $a$, under a smallness condition on the size of the shock \cite{MR4288384}. These methods have been extended to systems. In particular, the $a$-contraction was proven for sufficiently small shocks of the barotropic Navier-Stokes system \cite{MR4195742}. Using this strategy as a building block, a recent breakthrough allowed the proof of the convergence of solutions to the Navier-Stokes system to solutions to the Euler system in the framework of small-BV solutions in the vanishing viscosity limit \cite{NSvS}. The $a$-contraction theory also was a key tool in solving the long-standing problem of asymptotic stability of composite waves \cite{MR4560999}.

\subsubsection{Large shock waves} In the context of the viscous model, outside of the result of Kang-Vasseur on Burgers' equation, $a$-contraction results in the large shock setting (even in the scalar case) remain elusive. However, there are results for nonlinear stability of the viscous wave. Asymptotic orbital stability has been shown for a general scalar viscous conservation law \cite{MR0120469}. This was extended to the system setting (with full viscosity) by Sattinger \cite{MR0978372} \& Liu \cite{MR0791863}. However, in the system setting with partial viscosity (such as in the Navier-Stokes system), there are fewer results. In recent times, a  pointwise Green's function approach has been developed by Zumbrun and coauthors. This was carried out in the $1-$D scalar setting by Howard \cite{MR1699970}. Zumbrun and coauthors have worked on the development of conditions for whether these results extend to multidimensional systems with real viscosity, which may be found in \cite{MR2099037}.

\subsubsection{Related work} The closely related recent work of Kang shows $a$-contraction for some families of fluxes we address in this article (more specifically, he considers fluxes which are strictly convex and bounded by an exponential function $|A(x)| \leq Ce^{C|x|}$), but with smallness conditions on the shock size and the use of a different entropy function $\eta$ \cite{MR4188324}.  During the preparation of this paper, we learned about a new work on the $a$-contraction, allowing large shocks in the context of a cubic flux \cite{Wangcubic}. In this work, the failure of the convexity of $f$ is of importance. There, even for a single shock, the authors have to introduce an $a$ which is $C^2$. This suggests the possibility to stabilize in the sense of $a$-contraction shocks associated to more general fluxes. Our work shows that it is not possible to extend this result with a similar weight to arbitrary convex fluxes. It is of interest to understand what happens for other polynomial fluxes. Finally, we also mention the recent work of Kang-Oh showing $L^2$-contraction for the multidimensional Burgers' equation \cite{ContracMultD}.

\subsection{Results}
\par In this paper, we will investigate whether $a-$contraction holds for large shocks of viscous conservation laws. In particular, the idea is to understand cases in which there is stability of the viscous shock wave but the $a$-contraction fails. As a first result, in the scalar case, we provide a class of fluxes such that $a-$contraction fails to hold uniformly in the shock size for a large class of $a$, including the ones used by Stokols in \cite{MR4288384} and by Kang and Vasseur in \cite{MR4195742}. In particular, we will show that under some growth assumptions on $A'$, the shock connecting $0$ to $-K$ cannot be stabilized in the sense of $a$-contraction uniformly with respect to $K$.

\par This has two ramifications: firstly, when $a=1$, we note that this may be compared to the aforementioned result of Leger to showcase a ``viscous destabilization'' effect, meaning that the result in the inviscid case fails for the viscous case. Secondly, it shows that $a$-contraction can fail in situations where the viscous wave is asymptotically stable. This shows that the $a$-contraction property is stronger than nonlinear stability. 

\par Our proof is based on the same strategy as the negative result of Kang and Vasseur in \cite{MR3592682}. However, there are some crucial differences. We work with a larger class of weight functions $a$ rather than $a=1$. Further, we first fix a flux and then vary the shock size, while in their article, they fix a shock size then create a flux based on that. In the following, the shock connecting $0$ to $-K$ will be denoted $S$ (dropping the dependence on $K$) and the same will be done for its speed $\sigma$.
\begin{theorem}\label{maintheorem}
Let $A:\mathbb{R} \to \mathbb{R}$ be smooth, strictly convex, such that $A'(0)=0$ and such that for every real valued polynomial $P$ there exists $x$ in $\R_-$ such that $P(x)>A'(x)$. Fix $\tilde{a} \in W^{2,\infty}([0,1])$. Let $u_-=0$. There exists a sequence of positive numbers $(K_n)_{n \in \N}$ going to $\infty$ such that if there exists $n$ in $\N$ such that $u_+=-K_n$, then there exists some smooth initial data $u_0$ to (\ref{viscous}), with $u_0-S \in C_0^\infty$, such that for any Lipschitz shift $X$, there exists $T^*=T^*(\tilde{a}, \|\dot{X}\|_\infty) > 0$ such that for every $t$ in $(0,T^*)$ and with $a(x):=\tilde{a}(\frac{u_--S(x)}{u_--u_+})$:
\begin{equation*}
\int_\mathbb{R}a(x)|u(x+X(t),t)-S(x)|^2\diff x > \int_\mathbb{R} a(x)|u_0(x)-S(x)|^2\diff x \ ,
\end{equation*}
where $u$ is the solution to (\ref{viscous}) with initial data $u_0$.
\end{theorem}
Without loss of generality, we will assume throughout this paper that $A(0)=0$. Furthermore, $u_-$ could be any non-positive number instead of 0. \\
As a second result, we consider the barotropic Navier-Stokes system:
\begin{equation}\label{bns}
\begin{cases}
v_t-u_x=0, & (x,t) \in \mathbb{R} \times (0,\infty), \\ 
u_t+p(v)_x=\left(\frac{\mu(v)}{v}u_x\right)_x, & (x,t) \in \mathbb{R} \times (0,\infty), \\
(v(x,0),u(x,0))=(v_0(x), u_0(x)),  & x \in \mathbb{R},
\end{cases}
\end{equation}
where $v$ denotes the specific volume, $u$ is the velocity, and $p(v)$ is the pressure law. Here, we consider $p(v)=v^{-\gamma}$, $\gamma >1$, and the pressure dependent viscosity $\mu(v)=\gamma v^{-\gamma}$. 
This system also admits traveling wave solutions; given a speed $\sigma$ and states $(v_-,u_-)$ and $(v_+,u_+)$ that satisfy the Rankine-Hugoniot and Lax conditions:
\begin{equation}\label{rh}
\begin{cases}
    -\sigma(v_+-v_-)-(u_+-u_-)=0, \\
    -\sigma(u_+-u_-)+(p(v_+)-p(v_-))=0,
\end{cases}
\text{and } u_->u_+, \text{ and } \text{sgn}(\sigma)=\text{sgn}(v_+-v_-)
\end{equation}
there exists a viscous shock $(\tilde{v},\tilde{u})(x-\sigma t)$  satisfying
\begin{equation}\label{viscousshocksystem}
\begin{cases}
-\sigma \tilde{v}'-\tilde{u}'=0, \\
-\sigma \tilde{u}'+p(\tilde{v})'=\left(\frac{\mu(\tilde{v})}{v}\tilde{u}'\right)', \\
\lim_{\xi \to \pm \infty}=(\tilde{v}, \tilde{u})=(v_\pm, u_\pm), \\
\lim_{\xi \to \pm \infty}(\tilde{v}', \tilde{u}')=(0,0), \\
\sigma=-\sqrt{-\frac{p(v_+)-p(v_-)}{v_+-v_-}} \text{ if } v_->v_+, \\
\sigma=\sqrt{-\frac{p(v_+)-p(v_-)}{v_+-v_-}} \text{ if } v_-<v_+. 
\end{cases}
\end{equation}
As mentioned above, Kang and Vasseur showed $a$-contraction for small shock profiles using a weight parameterized as $a(x)=\tilde{a}\left(\frac{p(v_-)-p(\tilde{v})}{p(v_-)-p(v_+)}\right)$, for some $\tilde{a} \in W^{2,\infty}([0,1])$ \cite{MR4195742}. For convenience, we will denote:
\begin{align}
    b(\tilde{v}):=\frac{p(\tilde{v})-p(v_-)}{p(v_+)-p(v_-)},
\end{align}
so that the weight is parameterized as $a(x)=\tilde{a}(b(\tilde{v}(x)))$. Here, the contraction is obtained with respect to the functional associated with the BD-entropy (see \cite{MR1989675}, \cite{MR2257849}, \cite{MR1978317}):
\[
E((v_1,u_1),(v_2,u_2)):=\frac{1}{2}(u_1+p(v_1)_x-u_2-p(v_2)_x)^2+Q(v_1|v_2),
\]
where $Q(v_1|v_2):=Q(v_1)-Q(v_2)-Q'(v_2)(v_1-v_2)$ is the relative functional of the entropy flux $Q$ defined by $Q(v)=\frac{v^{-\gamma+1}}{-\gamma+1}$. Utilizing this, the authors were able to show uniqueness of entropy shocks for the isentropic Euler system in the class of inviscid limits of solutions to the Navier-Stokes system, again under a smallness condition for the size of the shock \cite{MR4228501}. As a second result, we show that $a$-contraction using a weight of this type fails for sufficiently large 1-shocks, i.e. when $v_+$ is sufficiently close to zero:
\begin{theorem}\label{maintheorem2}
Consider the system \eqref{bns}. Fix $\tilde{a} \in W^{2,\infty}([0,1])$. Let $(v_-,u_-) \in \mathbb{R}^+ \times \mathbb{R}$ be a given constant state. There exists $\epsilon > 0$ such that the following is true. Let $(v_+,u_+) \in \mathbb{R}^+ \times \mathbb{R}$ satisfy \eqref{rh}, and let $|v_+|<\eps$. Then for every $\delta>0$ and $s>0$ there exists smooth initial data $(v_0,u_0)$ to \eqref{bns} with $\tilde{U}-(v_0,u_0) \in [C_0^\infty(\mathbb{R})]^2$ and $\|\tilde{U}-(v_0,u_0)\|_{H^s}<\delta$ such that for any Lipschitz shift $X$ and every $t$ in $(0,T^*)$:
\[
\int_\mathbb{R}a(x)E(u(x+X(t),t)|\tilde{U}(x))\diff x > \int_\mathbb{R}a(x)E(u_0(x)|\tilde{U}(x))\diff x,
\]
for some $T^*=T^*(\tilde{a},\|\dot{X}\|_\infty)$, where $u$ is the solution to \eqref{bns} with initial data $(v_0,u_0)$. Here, $\tilde{U}=(\tilde{v}, \tilde{u})$ is the viscous shock connecting $(v_-,u_-)$ and $(v_+,u_+)$, and $a(x)=\tilde{a}(b(\tilde{v}(x)))$.
\end{theorem}
\subsection{Proof outlines}
For both theorems, the scheme of the proof is as follows.
First, we write the time derivative of the relative entropy:
\begin{align*}
    \frac{\diff}{dt}\int_\mathbb{R}\eta(u^{X}|S)\diff x,
\end{align*}
as $\dot{X}(t)(L(u^X-S)-Q_1(u^X-S))+Q_2(u^X-S)+R(u^X-S)$ where $Q_1$ and $Q_2$ are quadratic functionals on $H^1(\R)$, $L$ is linear, and $R$ is a nonlinear, higher order, remainder. We show that a weaker topology than $H^1$ makes $L$, $Q_1$, $Q_2$ and $R$ continuous, allowing us to consider at first a broader topological space. Crucially, the topology is coarse enough to allow for some functions equal to a nonzero constant at $\pm \infty$ to belong to the closure of test functions for that topology. In the scalar case, the norm is given by
\begin{align*}
    f \mapsto \sqrt{\int_{\R}|f(x)|^2|S'(x)|\diff x+\int_{\R}|f'(x)|^2\diff x} \, .
\end{align*}
Then, we show that there exists a solution $v$ in that space to the equation $L(v)=0$ which also satisfies that $Q_2(v)$ is positive. Finally, we apply a perturbative argument to conclude that for the original problem, there exists a solution $u$ such that 
\begin{align*}
    &L(u_0-S)=Q_1(u_0-S), \\
    &Q_2(u_0-S)+R(u_0-S)>0 \, ,
\end{align*}
and initially, $u_0-S$ is compactly supported. Even though both problem are solved through similar strategies, we will emphasize the differences that occur when going to the system case. In particular, in the system case, we have more freedom when choosing the perturbation.
\section{Preliminaries for scalar case}
Fix a smooth, strictly convex function $A : \R \rightarrow \R$ such that $A'(0)=0$ and such that for every polynomial function $P$, there exists $x$ in $\R_-$ such that $|A'(x)| \geq |P(x)|$. Fix $\tilde{a}$ in $W^{2,\infty}([0,1])$ such that $\tilde{a}$ is valued in $\R^+$.
\subsection{Existence and uniqueness of perturbations}
For \eqref{viscous}, it is well-known that one has global well-posedness with initial data $u_0 \in L^\infty$. For our purposes, we need finer bounds on the solution, starting with initial data $u_0-S \in C_0^\infty$, which will be proved in Section \ref{shiftapprox}.

\subsection{Relative entropy structure}\label{relentropyscalar}
In this part, we use the relative entropy method to give the desired representation for the time derivative. Firstly, we recall the following. Given a strictly convex entropy $\eta$ for \eqref{viscous}, we define the relative entropy by:
\begin{equation*}
    \eta(u|v)=\eta(u)-\eta(v)-\eta'(v)(u-v).
\end{equation*}
Also, we define the entropy flux $G$ of $\eta$:
\begin{equation*}
    G'=\eta'A',
\end{equation*}
the flux of the relative entropy:
\[
F(u;v)=G(u)-G(v)-\eta'(v)(A(u)-A(v)),
\]
and the relative flux:
\[
A(u|v)=A(v)-A(v)-A'(v)(u-v).
\]
Finally, for convenience, given any function $f(x,t)$, we will use the notation for the shifted function:
\begin{equation*}
    f^{X}(x,t):=f(x+X(t), t).
\end{equation*}
It is easy to check that the shifted viscous profile $S^X$ satisfies:
\[
S^X_t-\dot{X}(t)S^X_x+(A(S^X)-\sigma S^X)_x=S^X_{xx}.
\]
\begin{remark}
In \cite{MR3592682}, a very natural change of variables to the "shock variable" was introduced. For any functions $B,C$ such that the following integral is finite
\[
\int_{\R}C(S(x))B(u^{-X}(x,t)-S(x))S'(x)\diff x,
\]
where $S$ is the viscous profile satisfying \eqref{viscousshock}, $u$ is a solution to \eqref{viscous}, and $X:\mathbb{R}_{\geq 0} \to \mathbb{R}$ is some shift function, we may make a change of variables $y=S(x)$ (as $S$ is a monotone function) and define:
\begin{equation*}
    w(y,t)=u^{-X}(S^{-1}(y),t)-y,
\end{equation*}
Then, the integral becomes:
\[
\int_{u_-}^{u_+}C(y)B(w(y,t))\diff y.
\]
Working in the shock variable allows, among other things, the use of a nonlinear Poincar\'e inequality when proving contraction results. In particular, here, we will use it to make computations easier.
\end{remark}
With this in mind, we have the following result on the relative entropy. From here on out, we will fix the entropy $\eta(u)=u^2$. In the following, we have hidden the dependence of the functions $u$ and $w$ on $t$.
\begin{lemma}\label{structurelemma}
Let $S$ be the viscous shock solution defined by \eqref{viscousshock}, with $u_-=0, u_+=-K$. Let $u$ be a solution to \eqref{viscous} with initial data $u_0$ satisfying $u_0-S \in C^{\infty}_0(\R)$. Let $a(x)=\tilde{a}\left(\frac{-S(x)}{K}\right)$ for some $\tilde{a} \in W^{2,\infty}([0,1])$. Let $X:\mathbb{R}_{\geq 0} \to \mathbb{R}$ be $C^1$. Then, we have:
\begin{align}\label{representation}
\dv{}{t}\int_\mathbb{R}a\eta(u^{-X}|S)\diff x=&2\dot{X}(t)Y(w)+Z(w)+R_1(w),
\end{align}
where 
\begin{align*}
Y(w)=\int_{-K}^0\tilde{a}\left(\frac{-y}{K}\right)w(y)w_y(y)\diff y+\int_{-K}^0\tilde{a}\left(\frac{-y}{K}\right)w(y)\diff y.
\end{align*}
Further,
\begin{align*}
    Z(w)=&2\sigma\int_{-K}^0\tilde{a}\left(\frac{-y}{K}\right)w(y)\diff y \notag \\
&+2\int_{-K}^0\tilde{a}\left(\frac{-y}{K}\right)|w_y(y)|^2(-\sigma y+A(y))\diff y \notag \\
&+2\int_{-K}^0\tilde{a}\left(\frac{-y}{K}\right)A(w(y)+y|y)\diff y \notag \\
&-\frac{1}{K^2}\int_{-K}^0\tilde{a}''\left(\frac{-y}{K}\right)w(y)^2(-\sigma y+A(y))\diff y \notag \\
&+\frac{1}{K}\int_{-K}^0\tilde{a}'\left(\frac{-y}{K}\right)w(y)^2(-\sigma+2A'(y))\diff y, 
\end{align*}
and 
\begin{align*}
R_1(w)=\frac{1}{2K}\int_{-K}^0\tilde{a}'(\frac{-y}{K})\left(\int_y^{w+y}(G^{(3)}(k)-2yA^{(3)}(k)(w-y-k)^2)\diff k\right)\diff y.
\end{align*}
\end{lemma}
\begin{remark}
Note that we have a bound $|R_1(w)| \leq C||w||_\infty^3$ for $C$ depending only on $A, K$, and $||\tilde{a}||_{W^{2,\infty}}$. 
\end{remark}
\begin{proof}
We calculate:
\begin{align}\label{step1}
\dv{}{t}\int_{\mathbb{R}}a^X\eta(u|S^X)\diff x=\dot{X}(t)\int_{\mathbb{R}}a'^X\eta(u|S^X)\diff x+\int_\mathbb{R}a^X\eta(u|S^X)_t\diff x.
\end{align}
Now, expanding the time derivative of the relative entropy:
\begin{align}\label{tderivrelentropy}
\eta(u|S^X)_t=&(\eta'(u)-\eta'(S^X))u_t-\eta''(S^X)(u-S^X)S^X_t\notag \\
=&(\eta'(u)-\eta'(S^X))(u_{xx}-(A(u))_x) \notag \\
&-\eta''(S^X)(u-S^X)(S^X_t+S^X_{xx}+\sigma S^X_x-(A(S^X))_x).
\end{align}
Substituting \eqref{tderivrelentropy} into \eqref{step1}, and then adding and subtracting 
\[
\int_\mathbb{R}a^X(\eta'(u)-\eta'(S^X))S^X_{xx}\diff x, 
\]
gives:
\begin{align}\label{step2}
\dv{}{t}\int_{\R}\eta(u|S^X) \diff x=&\dot{X}(t)\left(\int_\mathbb{R}a'^X\eta(u|S^X)\diff x-\int_\mathbb{R}a^X\eta''(S^X)(u-S^X)S^X_x\diff x\right) \notag \\
&+\int_\mathbb{R}a^X(\eta'(u)-\eta'(S^X))(u_{xx}-S^X_{xx})\diff x \notag \\
&+\int_\mathbb{R}a^X(\eta')(u|S^X)S^X_{xx}\diff x \notag \\
&-\sigma \int_\mathbb{R}a^X\eta''(S^X)(u-S^X)S^X_x\diff x \notag \\
&-\int_\mathbb{R}a^X(\eta'(u)-\eta'(S^X))(A(u))_x\diff x\notag \\
&+\int_\mathbb{R}a^X\eta''(S^X)(u-S^X)(A(S^X))_x\diff x.
\end{align}
Next, we calculate:
\begin{align*}
(F(u;S^X))_x=\eta'(u)A'(u)u_x-\eta''(S^X)S^X_x(A(u)-A(S^X))-\eta'(S^X)(A(u))_x.
\end{align*}
Adding and subtracting $\int_{\R}a^X(F(u;S^X))_x\diff x$ in \eqref{step2}, we obtain:
\[
\dv{}{t}\int_{\mathbb{R}}a^X\eta(u|S^X)\diff x=\sum_{k=1}^6I_k,
\]
where:
\begin{align*}
I_1&=\dot{X}(t)\left(\int_\mathbb{R}a'^X\eta(u|S^X)\diff x-\int_\mathbb{R}a^X\eta''(S^X)(u-S^X)S^X_x\diff x\right), \\
I_2&=-\sigma \int_\mathbb{R}a^X\eta''(S^X)(u-S^X)S^X_x\diff x, \\
I_3&=\int_\mathbb{R}a^X(\eta'(u)-\eta'(S^X))(u_{xx}-S^X_{xx})\diff x, \\
I_4&=-\int_\mathbb{R}a^X(F(u;S^X))_x\diff x, \\
I_5&=-\int_\mathbb{R}a^X(\eta'(S^X))_xA(u|S^X)\diff x, \\
I_6&=\int_\mathbb{R}a^X(\eta')(u|S^X)S^X_{xx}\diff x. \\
\end{align*}
For the shift term $I_1$, we integrate the first term by parts to obtain:
\[
I_1=-\dot{X}(t)\left(\int_\mathbb{R}a\eta(u|S^X)_x\diff x+\int_\mathbb{R}a\eta''(S^X)(u-S^X)S^X_x\diff x \right).
\]
Now, integrating by parts for $I_4$ and expanding:
\begin{align*}
I_4=&\int_\mathbb{R}a'^XF(u;S^X)\diff x \\ =&\int_\mathbb{R}a'^X\frac{\eta''(S^X)A'(S^X)}{2}(u-S^X)^2\diff x \\
&+\frac{1}{2}\int_{\R}a'^X\left(\int_{S^X}^u(G^{(3)}(k)-2S^XA^{(3)}(k))(u-k)^2 \diff k\right)\diff x.
\end{align*}
Now, specifying the quadratic entropy $\eta(u)=u^2$, we obtain for $I_3$:
\begin{align*}
I_3&=2\int_\mathbb{R}a^X(u-S^X)(u_{xx}-S^X_{xx})\diff x \\
&=-2\int_\mathbb{R}a^X|(u-S^X)_x|^2\diff x+\int_\mathbb{R}a''^X(u-S^X)^2\diff x.
\end{align*}
So, all of our simplifications (alongside specifying $\eta(u)=u^2$) yield:
\begin{align}\label{step3}
\dv{}{t}\int_{\mathbb{R}}a^X\eta(u|S^X)\diff x =&-2\dot{X}(t)\left(\int_\mathbb{R}a^X(u-S^X)(u-S^X)_x\diff x+\int_\mathbb{R}a^X(u-S^X)S^X_x\diff x\right) \notag \\
&-2\sigma \int_\mathbb{R}a^X(u-S^X)S^X_x\diff x +\int_\mathbb{R}a''^{X}(u-S^X)^2\diff x \notag \\
&-2\int_\mathbb{R}a^X|(u-S^X)_x|^2\diff x +\int_\mathbb{R}a'^XA'(S^X)(u-S^X)^2 \diff x \notag \\
&-2\int_\mathbb{R}a^XA(u|S^X)S^X_x\diff x \notag \\
&+\frac{1}{2}\int_{\R}a'^X\left(\int_{S^X}^u(G^{(3)}(k)-2S^XA^{(3)}(k))(u-k)^2 \diff k\right)\diff x.
\end{align}
Making the change of variables $x \mapsto x-X(t)$, and then $y=S(x)$ in \eqref{step3}, we obtain exactly \eqref{representation}.
\end{proof}
\subsection{Reduction of the shift condition}
Now, in order to ensure that the time derivative is positive at $0$ for any shift $X$, we must construct an initial perturbation $\overline{w}:=u_0-S$ such that the term $Y(\overline{w})$ in Lemma \ref{structurelemma} equals zero. The next lemma shows that constructing an initial perturbation $\overline{w}$ that fits $Y(\overline{w})=0$ can be done by a perturbative argument from a function w satisfying a certain linear condition. The strategy is to slightly perturb $w$ by a function $\phi \in C_0^\infty([-K,0])$ that satisfies $\phi(K)=\phi(0)=0$.
\begin{lemma}\label{ift}
Let $w \in C_0^\infty([-K,0])$ such that $\int_{-K}^0\tilde{a}\left(\frac{-y}{K}\right)w(y)\diff y=0$. Let $\phi \in C_0^\infty([-K,0])$ such that $\int_{-K}^0\tilde{a}\left(\frac{-y}{K}\right)\phi(y) \diff y \neq 0$. \\
Let $\tilde{a}\left(\frac{-y}{K}\right)=\overline{a}(y)$. \\
With $\lambda_*=\frac{-\int_{-K}^0\frac{\overline{a}'}{2}w^2\diff y}{\int_{-K}^0\overline{a}\phi \diff y}$, there exists $\epsilon_0,\delta > 0$ and $\lambda:(0,\epsilon_0) \rightarrow (\lambda_*-\delta,\lambda_*+\delta)$ such that for every $\epsilon \in (0,\epsilon_0)$
\[
\int_{-K}^0\tilde{a}\left(\frac{-y}{K}\right)\overline{w}(y)\overline{w}_y(y)\diff y+\int_{-K}^0\tilde{a}\left(\frac{-y}{K}\right)\overline{w}(y)\diff y=0,
\]
where $\overline{w}=\epsilon w-\epsilon^2\lambda(\epsilon) \phi$, \\
and such that
\begin{equation}\label{lambdacont}
    \lim_{\epsilon \rightarrow 0} \lambda(\epsilon)=\lambda_* \, .
\end{equation}
\end{lemma}
\begin{proof}
We calculate:
\[
\int_{-K}^0\overline{a}\overline{w}\overline{w}_y\diff y=-\int_{-K}^0\frac{\overline{a}'}{2}(\epsilon^2w^2-2\epsilon^3\lambda \phi w+\epsilon^4\lambda^2\phi^2)\diff y.
\]
So, in order for this to be equal to $-\int_{-K}^0\overline{a}(\epsilon w-\epsilon^2\lambda \phi) \diff y$, it suffices to find $\lambda$ such that the function:
\[
f(\epsilon, \lambda)=\lambda \int_{-K}^0\overline{a}\phi \diff y+\int_{-K}^0\frac{\overline{a}'}{2}w^2\diff y-\lambda \epsilon \int_{-K}^0\overline{a}'\phi w \diff y+\lambda^2\epsilon^2\int_{-K}^0\frac{\overline{a}'}{2}\phi^2\diff y,
\]
is equal to zero. We do so via the implicit function theorem. Firstly, note that for $\epsilon=0$, we may choose $\lambda_*=\frac{-\int_{-K}^0\frac{\overline{a}'}{2}w^2\diff y}{\int_{-K}^0\overline{a}\phi \diff y}$. Secondly, we verify:
\[
D_\lambda f(0, \lambda_*)=\int_{-K}^0\overline{a}\phi \diff y \neq 0.
\]
So, by the implicit function theorem, we obtain the result.
\end{proof}
\begin{remark}
Note that by \eqref{lambdacont}, we may choose $\epsilon$ sufficiently small such that the $\lambda$ given by Lemma \ref{ift} is arbitrarily close to $\lambda_*$. More specifically, for $\eps$ sufficiently small, $\lambda(\epsilon)=\lambda_*+\beta(\epsilon)$, where $\beta(\epsilon) \to 0$ as $\epsilon \to 0$. 
\end{remark}
\subsection{Estimates for $A$}
In this section, we state and prove the estimates needed for the flux function $A$. 
\begin{lemma}\label{unbounded}
Fix $\theta \in (0,1)$. There exists some $M_{\theta} > 0$ such that the function:
\[
\phi_{\theta}:=\frac{A'(\cdot)}{A'(\theta \cdot)},
\]
is well-defined on some interval $(-\infty,-M_\theta]$ and is unbounded on $(-\infty,-M_\theta]$.
\end{lemma}
\begin{proof}
As $A'$ is increasing and that for every $M>0$ there exists $x$ in $\R_-$ such that $-M>A'(x)$, we have that $A'(x)$ goes to $-\infty$ as $x$ goes to $-\infty$. \\
Hence, $A'(x)$ goes to $-\infty$ as $x$ goes to $\infty$, so there exists $M_\theta$ such that $A'(\theta x)$ is negative for $x \leq -M_\theta$. By contradiction, assume that $\phi_{\theta}$ is bounded on $(-\infty,-M_\theta]$. Then there exists $C$ positive such that for any $x \in (-\infty, -M_{\theta}]$
\begin{align*}
\ln(|A'(x)|)-\ln(|A'(\theta x)|) \leq C \, .
\end{align*}
Thus, with $C'=\sup \lbrace |A'(x)| \, : \, x \in [-\theta^{-1}M_{\theta},-M_{\theta}] \rbrace$, we have that for any $n$ in $\N \cup \lbrace 0 \rbrace$ and $x$ in $[-\theta^{-(n+1)}M_{\theta},-\theta^{-n}M_{\theta}]$
\[
\ln(|A'(x)|) \leq C'+nC \leq C'+\frac{(\ln(|x|)-\ln(M_{\theta}))C}{|\ln(\theta)|} \, .
\]
Hence, for every $x$ in $(-\infty,-M_\theta]$, $A'(x) \geq -e^{C'+nC}$ where $n=\frac{\ln(|x|)-\ln(M_{\theta})}{|\ln(\theta)|}$. \\
Denote $p$, some integer greater than $\frac{C}{|\ln(\theta)|}$. \\
Thus, for every $x$ in $(-\infty,-M_\theta]$, $A'(x) \geq -e^{C'}e^{\frac{C\ln(M_{\theta})}{\ln(\theta)}}|x|^p$. A contradiction. \\
Thus, $\phi_{\theta}$ is unbounded on $(-\infty,-M_\theta]$.
\end{proof}
\begin{remark}
We will use $\theta$ very close to $1$. Note that the lemma for any $\theta_1$ implies the lemma for any $0<\theta_2 \leq \theta_1$ on the interval $(-\infty,-M_{\theta_2}]$.
\end{remark}
One consequence of this is that $\sigma$ is much smaller than $A'$ for $u_+$ near $-\infty$. \\
Recall that $\sigma=\frac{A(0)-A(-K)}{K}$.
\begin{lemma}\label{boundsigmaaa'}
We have that for any $\theta \in (0,1)$, $\rho$ and $K_0$ positive there exists $K \geq K_0$ such that
\begin{align*}
    \rho A'(-K) \leq A'(-\theta K) < 0 \, ,
\end{align*}
and
\begin{align}\label{sigmabound}
\left|\frac{\sigma}{A'(-K)}\right|\leq \rho \, , 
\end{align}
and such that for every $y$ in $\left[-\frac{K}{2},0\right]$
\begin{align}\label{a'bound}
    |A'(y)| \leq \rho |A'(-K)| \, ,
\end{align}
and such that for every $x$ in $[-K,0]$
\begin{align}\label{abound}
    |A(z)| \leq \rho K |A'(-K)| \, .
\end{align}
\end{lemma}
\begin{proof}
Fix $K_0$ positive, $\rho$ in $\left(0,\frac{1}{2}\right]$ and $\theta$ in $(0,1)$. \\
For any $\theta'$ in $(0,1)$
\begin{align*}
    |\sigma|=&-\frac{1}{K}\int_{-K}^0A'(y)\diff y=-\frac{1}{K}\left(\int_{-K}^{-\theta' K}A'(y) \diff y+\int_{-\theta' K}^0A'(y) \diff y\right) \\ \leq& \frac{1}{K}(|A'(-K)|(1-\theta')K+|A'(-\theta' K)|\theta' K),
\end{align*}
Now, using Lemma \ref{unbounded}, for any $\theta'$ in $\left( 0,\frac{1}{2} \right)$, for any $K_0$ positive, there exists $K \geq K_0$ such that  $|A'(-K)| \geq \frac{1}{\rho} |A'\left( -\theta' K \right)|$. \\
Hence, there exists $K \geq K_0$ such that $|A'(-K)| \geq \frac{1}{\rho}|A'\left( -\max(\theta',\theta) K \right)|$. Thus, for such a $K$, we have
\begin{align*}
    \frac{1}{K}(|A'(-K)|(1-\theta')K+|A'(-\theta' K)|\theta' K) \leq& \frac{1}{K}(|A'(-K)|(1-\theta')K+\rho \theta' K |A'(-K)|)  \\ \leq& |A'(-K)|(\rho \theta' +(1-\theta')) \, ,
\end{align*}
and
\begin{align*}
    \rho|A'(-K)| \geq |A'(-\theta K)| \, ,
\end{align*}
and, as $A'$ is decreasing on $\R_-$
\begin{align*}
    |A'(-K)| \geq \frac{1}{\rho} \left|A'\left(-\frac{K}{2}\right)\right| \, .
\end{align*}
Thus, for any $K_0$ positive, there exists $K \geq K_0$ such that
\begin{align*}
    |\sigma| \leq 2\rho |A'(-K)| \, ,
\end{align*}
and so, as $A$ is decreasing on $\R_-$ and as $|A(-K)|=|\sigma K|$, we get the result. 
\end{proof}
\section{Proof of Theorem \ref{maintheorem}}
In this section, we prove Theorem \ref{maintheorem}. The strategy is to make a small initial perturbation that causes the derivative of the relative entropy given by the representation \eqref{representation} to be positive at $t=0$. First, fix $K_0$ positive. Fix $K \geq K_0$ such that Lemma \ref{boundsigmaaa'} applies for $\rho=\min\left(\frac{9}{10},\frac{\inf\tilde{a}}{48\sup\tilde{a}},\frac{3}{200(1+\|\tilde{a}\|_{\infty}+\|\tilde{a}'\|_{\infty})}\right)$ and for some $\theta$ in $\left(\frac{1}{2},1\right)$ such that $\min\tilde{a} \geq \frac{9}{10}\tilde{a}(1)$.
\subsection{Small perturbation of $S$}\label{smallperturbation} Firstly, assume that we can construct a perturbation $w$ with $\int_{-K}^0\tilde{a}(\frac{-y}{K})w(y)\diff y=0$. Then, for $\epsilon < \epsilon_0$, consider the perturbation $\overline{w}$ constructed in Lemma \ref{ift}, where we also take $\phi$ in the lemma such that $|\phi|,|\phi'| \leq C$ independently of $K$. Doing the Taylor expansion of the relative flux $A(\overline{w}(y)+y|y)=A(\epsilon w(y)+\epsilon^2\lambda\phi(y)+y|y)$, we obtain for every $y$ in $(-K,0)$:
\[
A(\overline{w}(y)+y|y)=\frac{1}{2}A''(y)\overline{w}(y)^2+\frac{1}{2}\int_y^{\overline{w}+y}A^{(3)}(k)(\overline{w}-y-k)^2\diff k.
\]
Consider the solution $u$ to \eqref{viscous} with initial data $S+\overline{w}$. Using the representation \eqref{representation} with $u$ and for any $X$ in $C^1([0,\infty),\R)$, at $t=0$ we have
\begin{align*}
\dv{}{t}\int_\mathbb{R}a\eta(u^{-X}|S)^2\diff x|_{t=0}=&2\sigma\int_{-K}^0\tilde{a}\left(\frac{-y}{K}\right)\overline{w}(y)\diff y \notag \\
&+2\int_{-K}^0\tilde{a}\left(\frac{-y}{K}\right)|\overline{w}_y(y)|^2(-\sigma y+A(y))\diff y \notag \\
&+\int_{-K}^0\tilde{a}\left(\frac{-y}{K}\right)A''(y)\overline{w}(y)^2 \diff y\notag \\
&-\frac{1}{K^2}\int_{-K}^0\tilde{a}''\left(\frac{-y}{K}\right)\overline{w}(y)^2(-\sigma y+A(y))\diff y \notag \\
&+\frac{1}{K}\int_{-K}^0\tilde{a}'\left(\frac{-y}{K}\right)\overline{w}(y)^2(-\sigma+2A'(y))\diff y \notag \\
&+\int_{-K}^0\tilde{a}\left(\frac{-y}{K}\right)\left(\int_y^{\overline{w}+y}A^{(3)}(k)(\overline{w}+y-k)^2\diff k\right)\diff y \notag \\
&+R_1(w) \notag.
\end{align*}
For the term linear in $\overline{w}$, recalling the definition of $\overline{w}$ from Lemma \ref{ift} and the remark thereafter, we see that:
\begin{align*}
2\sigma\int_{-K}^0\tilde{a}\left(\frac{-y}{K}\right)\overline{w}(y) \diff y=&-2 \sigma \epsilon^2(\lambda_*+\beta(\epsilon))\int_{-K}^0\tilde{a}\left(\frac{-y}{K}\right)\phi \diff y \\
=&-\frac{\epsilon^2\sigma}{K} \int_{-K}^0\tilde{a}'\left(\frac{-y}{K}\right)w(y)^2\diff y \\
&-2 \sigma \epsilon^2\beta(\epsilon)\int_{-K}^0\tilde{a}\left(\frac{-y}{K}\right)\phi(y) \diff y.
\end{align*}
So, expanding the powers of $\overline{w}$ for the other terms, and then pulling out $\epsilon^2$, we obtain, for a fixed $K$:
\begin{align}\label{eperturbation}
\dv{}{t}\int_\mathbb{R}a|u^{-X}-S|^2\diff x|_{t=0}=&\epsilon^2\biggl(\int_{-K}^0\tilde{a}\left(\frac{-y}{K}\right)A''(y)w(y)^2\diff y \notag \\
&+\frac{2}{K}\int_{-K}^0\tilde{a}'\left(\frac{-y}{K}\right)(-\sigma+A'(y))w(y)^2\diff y \notag \\
&\;-\frac{1}{K^2}\int_{-K}^0\tilde{a}''\left(\frac{-y}{K}\right)(-\sigma y+A(y))w(y)^2\diff y \\
&\; +2\int_{-K}^0\tilde{a}\left(\frac{-y}{K}\right)(-\sigma y+A(y))|w_y(y)|^2\diff y \notag \\
&+\frac{R_1(w)+R_2(w)+R_3(w)}{\eps^2}-2 \sigma \beta(\epsilon)\int_{-K}^0\tilde{a}\left(\frac{-y}{K}\right)\phi(y) \diff y \biggr)\notag ,
\end{align}
where:
\[
R_2(w)=\int_{-K}^0\tilde{a}\left(\frac{-y}{K}\right)\left(\int_y^{\overline{w}+y}A^{(3)}(k)(\overline{w}+y-k)^2\diff k\right)\diff y, 
\]
and:
\begin{align*}
R_3(w)=&2\int_{-K}^0\tilde{a}\left(\frac{-y}{K}\right)|\overline{w}_y(y)|^2(-\sigma y+A(y))\diff y \notag \\
&+\int_{-K}^0\tilde{a}\left(\frac{-y}{K}\right)(A''(y)\overline{w}(y)^2 \diff y\notag \\
&-\frac{1}{K^2}\int_{-K}^0\tilde{a}''\left(\frac{-y}{K}\right)\overline{w}(y)^2(-\sigma y+A(y))\diff y \notag \\
&+\frac{1}{K}\int_{-K}^0\tilde{a}'\left(\frac{-y}{K}\right)\overline{w}(y)^2(-\sigma+2A'(y))\diff y \notag \\
&-2\int_{-K}^0\tilde{a}\left(\frac{-y}{K}\right)|w_y(y)|^2(-\sigma y+A(y))\diff y \notag \\
&-\int_{-K}^0\tilde{a}\left(\frac{-y}{K}\right)A''(y)w(y)^2 \diff y\notag \\
&+\frac{1}{K^2}\int_{-K}^0\tilde{a}''\left(\frac{-y}{K}\right)w(y)^2(-\sigma y+A(y))\diff y \notag \\
&-\frac{1}{K}\int_{-K}^0\tilde{a}'\left(\frac{-y}{K}\right)w(y)^2(-\sigma+2A'(y))\diff y.
\end{align*}
Note that $|\frac{R_1(w)+R_2(w)+R_3(w)}{\eps^2}-2 \sigma \beta(\epsilon)\int_{-K}^0\tilde{a}(\frac{-y}{K})\phi(y) \diff y| \to 0$ as $\eps \to 0$. 
\subsubsection{Construction of $w$}
Define the following functional $F$:
\begin{align*}
F(w)&:=J_1(w)+J_2(w),
\end{align*}
where:
\begin{align*}
J_1(w):=&\int_{-K}^0\tilde{a}\left(\frac{-y}{K}\right)A''(y)w(y)^2\diff y+\frac{2}{K}\int_{-K}^0\tilde{a}'\left(\frac{-y}{K}\right)\left(-\sigma+A'(y)\right)w(y)^2\diff y \notag \\
&-\frac{1}{K^2}\int_{-K}^0\tilde{a}''\left(\frac{-y}{K}\right)\left(-\sigma y+A(y)\right)w(y)^2\diff y, \notag \\
J_2(w):=&2\int_{-K}^0\tilde{a}\left(\frac{-y}{K}\right)(-\sigma y+A(y))|w_y(y)|^2\diff y.
\end{align*}
Due to the expansion \eqref{eperturbation}, it suffices to show that for every $K_0$ there exists $K \geq K_0$ such that we may find a smooth function $w$ on $[-K,0]$ such that $F(w) > 0$. There exists a non-increasing smooth function $\tilde{w}:[-1,0] \to \mathbb{R}$ with the following properties:
\begin{enumerate}
    \item $\tilde{w}(x)=1$ for $x \in [-1,\frac{-1}{2}]$,
    \item $\tilde{w}(x)=-C$ for $x \in [-\frac{1}{4},0]$, where $C$ depends only on $\tilde{a}$,
    \item $\tilde{w}$ is nonincreasing on $[-1,0]$,
    \item $\int_{-1}^0\tilde{a}(-y)\tilde{w}(y)\diff y=0$.
\end{enumerate} 
\begin{remark}
There exists such a $\tilde{w}$ with $C$ between $\frac{\inf \tilde{a}}{\sup\tilde{a}}$ and $\frac{3\sup\tilde{a}}{\inf\tilde{a}}$. We choose $\tilde{w}$ such that $C$ is in that interval.
\end{remark} 
Then, let $w(y)=\tilde{w}\left(\frac{y}{K}\right)$. We will show that this function satisfies $F(w) > 0$ for $K \geq K_0$ chosen as above. Then, we will approximate $w$ by functions that are zero on the boundary in the topology induced by the functional $F$. \\
\subsubsection{Control on $J_1$}
We split the integral into two parts: one close to $0$, and one close to $-K$. Firstly, we compute:
\begin{align}\label{bottomhalf}
&\int_{-K}^{\frac{-K}{2}}\tilde{a}\left(\frac{-y}{K}\right)A''(y)w(y)^2\diff y+\frac{2}{K}\int_{-K}^{\frac{-K}{2}}\tilde{a}'\left(\frac{-y}{K}\right)(-\sigma+A'(y))w(y)^2\diff y \notag \\
&-\frac{1}{K^2}\int_{-K}^{\frac{-K}{2}}\tilde{a}''\left(\frac{-y}{K}\right)(-\sigma y+A(y))w(y)^2\diff y \notag \\ 
=&-\left(A'\left(\frac{-K}{2}\right)\right)\tilde{a}\left(\frac{1}{2}\right)+(A'(-K))\tilde{a}(1) \\ 
&+\frac{1}{K}\left(A\left(\frac{-K}{2}\right)+\frac{\sigma 
 K}{2}\right)\tilde{a}'\left(\frac{1}{2}\right)-\frac{1}{K}(A(-K)+ 
\sigma K)\tilde{a}'(1)-\frac{2\sigma}{K}\left(\tilde{a}\left(\frac{1}{2}\right)-\tilde{a}(1)\right) \notag \\
&+2\int_{-K}^{\frac{-K}{2}}\tilde{a}\left(\frac{-y}{K}\right)A''(y)\diff y. \notag 
\end{align}
Further, with $\theta$ as above
\begin{align}
&\int_{-K}^{\frac{-K}{2}}\tilde{a}\left(\frac{-y}{K}\right)A''(y)\diff y \geq \int_{-K}^{-\theta K}\tilde{a}\left(\frac{-y}{K}\right)A''(y)\diff y \notag  
 \\
&\geq \left(A'(-\theta K)-A'(-K)\right)\min_{[\theta,1]}\tilde{a} \, .\notag  
\end{align}
Thus, using Lemma \ref{boundsigmaaa'}
\begin{align}\label{goodpart}
&\int_{-K}^{\frac{-K}{2}}\tilde{a}\left(\frac{-y}{K}\right)A''(y)\diff y \geq \frac{9}{10}(\rho-1)\tilde{a}(1)A'(-K) =\frac{9}{10}(1-\rho)\tilde{a}(1)|A'(-K)|. 
\end{align}
In light of the bounds \eqref{sigmabound}, \eqref{a'bound}, \eqref{abound}, and \eqref{goodpart}, we find that \eqref{bottomhalf} is bounded from below by:
\[
\frac{1}{2}\tilde{a}(1)|A'(-K)| \, ,
\]
For the other part,using the bounds \eqref{sigmabound}, \eqref{a'bound}, \eqref{abound}, we compute:
\[
\frac{2}{K}\int_{\frac{-K}{2}}^{0}\tilde{a}'\left(\frac{-y}{K}\right)(-\sigma+A'(y))w(y)^2\diff y \geq -C\rho|A'(-K)|,
\]
and:
\[
-\frac{1}{K^2}\int_{\frac{-K}{2}}^0\tilde{a}''\left(\frac{-y}{K}\right)(-\sigma y+A(y))w(y)^2\diff y \geq -C\rho |A'(-K)|.
\]
So, in total, we see:
\begin{align}
J_1(w) \geq -\frac{1}{4}\tilde{a}(1)A'(-K) \, . \label{J1}
\end{align}
\subsubsection{Control on $J_2$}
For $J_2$, we recall that $w(y)=\tilde{w}\left(\frac{y}{K}\right)$, so we have:
\[
\int_{-K}^0\tilde{a}\left(\frac{-y}{K}\right)(-\sigma y+A(y))|w_y(y)|^2\diff y \geq \frac{-C}{K^2}\int_{-\frac{K}{2}}^{-\frac{K}{4}}|-\sigma y +A(y)|\diff y \geq -C\rho |A'(-K)| \, .
\]
This implies:
\begin{align}
J_2(w) \geq -C \rho |A'(-K)| \, . \label{J2}
\end{align}
Recall that $C$ depends only on $\tilde{a}$.
\subsubsection{Conclusion}
Hence, for every $K_0$ positive there exists $K \geq K_0$ and there exists $\epsilon_1$ in $(0,\epsilon_0)$ such that for every $\epsilon$ in $(0,\eps_1)$ we have that for every $X$ in $C^1([0,\infty),\R)$, the solution $u$ to \eqref{viscous} with initial data $S+\overline{w}$ satisfies that
\begin{align*}
    \dv{}{t}\int_\mathbb{R}a\eta(u^{-X}|S)^2\diff x|_{t=0} \geq \frac{\eps^2}{8} \tilde{a}(1)|A'(-K)|
\end{align*}
\subsection{Approximation in the functional topology}\label{initdataconstruct}
In light of \eqref{J1} and \eqref{J2}, we see that for every $K_0$ there exists $K \geq K_0$ such that there exists $\overline{w}$ such that $F(\overline{w}) >0$ and $Y(\overline{w})=0$. Now, we show that we can find for such $K$ a function $q \in C_0^\infty((-K,0))$, with $F(q)>0$, such that $Y(q)=0$. It suffices to approximate the perturbation $w$ by a function in $C_0^\infty$. Firstly, for $K$ fixed, consider the following norm
\[
\|f\|_H^2=\|f\|^2_{L^2([-K,0])}+\int_{-K}^0y(K+y)|f'(y)|^2\diff y,
\]
and define $H$ to be the closure of $C_0^\infty((-K,0))$ for this topology. $F$ is continuous on $H$ equipped with the induced topology. \\
Notice that $w$ is an element of $H$. In fact, one can first notice that $H^1_0(-K,0)$ is a subset of $H$ as $C_0^\infty((-K,0))$ is a dense subset of $H^1_0(-K,0)$ for the $H^1$ topology, which is stronger than the one induced by $\|\cdot\|_H$. Furthermore, fix $\delta > 0$ and define $\phi_{\delta}:[-K,0] \rightarrow \R$ as follows:
\[
\phi_{\delta}(x)=\begin{cases}
    (K+x)^{\delta} & x \in (-K, -K+\delta], \\
    \delta^{\delta} & x \in (-K+\delta,-\delta), \\
    x^\delta & x \in [-\delta, 0).
\end{cases}
\]
The family $(\phi_{\delta}w)_{\delta > 0}$ is a family of elements of $H^1_0$ converging to $w$ in $H$ when $\delta$ goes to 0. As $w$ is thus in $H$ and has integral $\int_{-K}^0\tilde{a}\left(\frac{-y}{K}\right)w(y)dy=0$, there must be a sequence of functions in $C^{\infty}_0(\R)$ all satisfying that their integral against $\tilde{a}\left(\frac{-y}{K}\right)$ is 0 converging to $w$ in $H$. In fact, it is a consequence of the following elementary lemma:
\begin{lemma}\label{kernellemma}
    Let $(X,\|\cdot\|_X)$ be a Banach space and let $f$ be a continuous linear form on $X$. If $D$ is a dense linear subspace of $X$, then $D \cap \ker(f)$ is a dense linear subspace of $\ker(f)$.
\end{lemma}
\begin{proof}
    Consider such a Banach space $(X,\|\cdot\|_X)$ and $f$ a continuous linear form on $X$. If $f$ is the null form then it is clear. If $f$ is not the null form, there exists $x \in D$ such that $f(x)=1$. Consider $y \in X$ such that $f(y)=0$. There exists $(d_n)_n \in D^{\mathbb{N}}$ such that $(d_n)_n$ converges towards $y$ in $X$ as $n$ goes to infinity. Hence, $(f(d_n))_n$ converges towards $0$ as $n$ goes to infinity. Hence, $(d_n-f(d_n)x)_n$ is a sequence valued in $D$ that converges towards $y$ as $n$ goes to infinity.
\end{proof}

\subsection{Approximation of Lipschitz shifts}\label{shiftapprox}
Finally, we need to show that Theorem \ref{maintheorem} holds for shifts that are merely Lipschitz and not $C^1$. We do so by approximation. Let $X \in W^{1,\infty}(\mathbb{R}_{\geq 0})$, $X(0)=0$. Without loss of generality (because we consider $T^*$ arbitrarily small in Theorem \ref{maintheorem}), assume that $X(t)=0$ for $t \geq 1$. Thus, $X \in W^{1,1}(\mathbb{R}_{\geq 0})$. Then, we have the following:
\begin{lemma}\label{c1shifts}
There exists a family $(X_n)_{n \in \mathbb{N}} \in C^1(\mathbb{R}_{\geq 0})^{\mathbb{N}}$ such that:
\begin{enumerate}
    \item $\|\dot{X_n}\|_\infty \leq \|\dot{X}\|_\infty$,
    \item $\|X_n-X\|_\infty \to 0$,
    \item $X_n(0)=0$ for all $n$.
\end{enumerate}
\end{lemma}
\begin{proof}
Extend $X$ to $\mathbb{R}$ by $0$. Then, define $X_n:=\phi_{\frac{1}{n}}*X$, where $(\phi_{\frac{1}{n}})_{n \in \mathbb{N}}$ is the following family of mollifiers with $\|\phi_{\frac{1}{n}}\|_1=1$:
\begin{equation*}
\phi_{\frac{1}{n}}(x)=n\eta(n(x+\frac{1}{n})),
\end{equation*}
where $\eta$ is the standard mollifier:
\begin{equation*}
\eta(x)=\begin{cases}
Ce^{\frac{1}{|x|^2-1}} & |x| < 1 \\
0 & |x| \geq 1,
\end{cases}
\end{equation*}
and $C$ is such that $\int_{\R}\eta(x)\diff x=1$. Note that convolution with the offset mollifier enforces the condition $X_n(0)=0$. Further, $\|X_n-X\|_{W^{1,1}(\mathbb{R})} \to 0$. The first claim follows by Young's convolution inequality, and the second by Sobolev embedding. Restrict $X_n$ to $\mathbb{R}_{\geq 0}$ to obtain the result. 
\end{proof}
Further, we need strong bounds on the solution in time, in order to show that the time derivative of the relative entropy is positive not just at $t=0$, but also for a short time afterwards. 
\begin{lemma}\label{bounds}
Let $u$ be the solution to \eqref{viscous} such that $u_0-S \in C_0^\infty(\mathbb{R})$. Then, there exists a positive constant $T$ such that $u-S \in C([0,T];H^2(\mathbb{R}))$. 
\end{lemma}
\begin{proof}
Firstly, we show that $u-S \in C([0,T];L^2(\mathbb{R}))$. We have that $u-S$ solves the equation:
\begin{equation}
\begin{cases}
(u-S)_t+(A(u))_x-(A(S))_x+\sigma S_x=(u-S)_{xx}, & (x,t) \in \mathbb{R} \times (0,\infty), \\ 
(u-S)(x,0)=u_0(x)-S(x), & x \in \mathbb{R}.
\end{cases}
\end{equation}
So, we have a representation using the Duhamel formula:
\[
u(\cdot,t)-S=e^{\partial^2_x t}(u_0-S)-\int_0^t\left(e^{\partial^2_x (t-s)}(A(u(s))-A(S))\right)_xds-\sigma \int_0^t e^{\partial^2_x (t-s)}S_xds.
\]
Using $L^1$ bounds on the heat kernel and its spatial derivative, alongside the maximum principle (so $u-S \in L^\infty$, with a uniform bound on the $L^\infty$ norm to handle the $A(u)-A(S)$ term) then applying Gronwall's inequality, we obtain $u(\cdot, t)-S \in L^2(\mathbb{R})$, with a bound (for some constants $C$ and $\omega$ depending on $u_0$):
\[
\|u(\cdot, t)-S\|_2 \leq C\left( e^{ \omega t}\|u_0-S\|_2 + 1 \right).
\]
It is then straightforward to check the strong continuity at $0$ using the Duhamel formula again. Next, we show that $u_x-S_x \in C([0,T];L^2(\mathbb{R}))$. As $S_x \in L^2(\mathbb{R})$, it suffices to work directly on the equation for $u_x$, and show $u_x \in C([0,T];L^2(\mathbb{R}))$. Differentiating \eqref{viscous}, we find that $v:=u_x$ satisfies the equation:
\begin{equation}\label{firstderiv}
\begin{cases}
v_t+(A'(u)v)_x=v_{xx}, & (x,t) \in \mathbb{R} \times (0,\infty), \\ 
v(x,0)= (u_{0})_{x}(x)=:v_0(x), & x \in \R.
\end{cases}
\end{equation}
Again, this gives a representation using the Duhamel formula:
\[
v(\cdot,t)=e^{\partial^2_x t}v_0-\int_0^t\left((e^{\partial^2_x (t-s)})(A'(u(s))v(s))\right)_xds.
\]
Using the $L^1$ bounds on the heat kernel and its spatial derivative (and also the maximum principle again, to uniformly bound $\|A'(u)\|_\infty$), alongside Gronwall again, we find that $v(\cdot, t) \in L^2(\mathbb{R})$, with a similar bound. Finally, we can check the strong continuity at $0$ using the Duhamel formula again.  \\
We obtain the uniform bound on the $L^2$ norm of the second derivative by differentiating \eqref{firstderiv} and applying the same techniques as above.
\end{proof}
Finally, we need a bound on the $L^2$ norms of various functions minus translates of themselves, depending on the $L^\infty$ norm of the translate. We will apply this to the functions $u(\cdot,t)$ and $u_x(\cdot,t)$ in the scalar case, and their analogs in the system case. Note that while $u(\cdot,t)$ is not in $L^2$, $u(\cdot,t)$ minus a translate of itself is.
\begin{lemma}\label{shifttranslate}
Let $w$ be in $H^1_{loc}(\R)$ and such that $w_x \in L^2(\R)$. Then, for every $m \in \mathbb{R}$
\[
\|w(\cdot)-w(\cdot+m)\|_2 \leq \|w_x\|_2|m|.
\]
\end{lemma}
\begin{proof}
First, notice that
\[
w(\cdot)-w(\cdot+m)=-\int^1_0w_x(\cdot+\theta m)m\diff \theta,
\]
and so, by Jensen's inequality:
\[
\int_{\R}|w(x)-w(x+m)|^2\diff x \leq m^2\int_{\R}\int^1_0|w_x(x+\theta m)|^2\diff \theta \diff x.
\]
We calculate using Fubini's Theorem:
\[
\|w(\cdot)-(w(\cdot+m)\|_2^2 \leq \|w_x\|_2^2m^2.
\] 
Taking the square root gives the result. 
\end{proof}
\begin{remark}
Note that we have uniform bounds on $\|u_x(\cdot,t)\|_2, \|u_{xx}(\cdot,t)\|_2$ for $0 \leq t \leq T^*$ for some fixed $T^*$ by Lemma \ref{bounds}. So Lemma \ref{shifttranslate} may be applied with $w(\cdot)=u(\cdot, t), w=u_x(\cdot,t)$.
\end{remark} 
Now, consider the family $\{X_n\}_{n=1}^\infty$ of shifts constructed in Lemma \ref{c1shifts}. For the initial perturbation constructed in Section \ref{initdataconstruct}, we have:
\[
\dv{}{t}\int_\mathbb{R} a(x)|u(x+X_n(t),t)-S(x)|^2\diff x|_{t=0}=\dot{X_n}(0)Y_n(0)+Z_n(0) \geq c > 0,
\]
for some $c$ uniform in $n$. Next, we use the following general lemma (which will also be used in the system case) to show that this time derivative remains positive uniformly in $n$ for a short time after $t=0$. For brevity, the proof is relegated to the appendix.
\begin{lemma}\label{approximationlemma}
Let $N=1,2$. Let $\tilde{w} \in \left( L^\infty(\R) \cap C^2(\R) \right)^N$, such that \newline $\tilde{w}' \in \left(L^2(\R) \cap L^\infty(\R)\right)^N$ and $\tilde{w}'' \in \left(L^2(\R)\right)^N$. Let $w:\R \times [0,T] \to \R^N$ such that $w-\tilde{w} \in C([0,T];\left(H^2(\R)\right)^N)$. Denote by $V$ and $\tilde{V}$ open sets containing respectively the range of $w$ and of $\tilde{w}$. Let $a(x)=\tilde{a}(\tilde{w}_1(x))$ for some $\tilde{a} \in W^{2,\infty}(\tilde{V})$. Let $q_0 \in C^1(\tilde{V},\R^N)$, $\tilde{q}_0 \in C^1(\tilde{V},\R)$, $q_1 \in C^1(V\times \tilde{V},\R^N)$ and $q_i \in C^1(V\times \tilde{V}, \R), 2 \leq i \leq 6$, verify the following properties:
\begin{enumerate}
\item There exists $q^5 \in C^2(V \cup \tilde{V})$ such that $q_5(u,v)=q^5(u_1)-q^5(v_1)$, and the Hessian of $q^5$ is uniformly bounded on $V \cup \tilde{V}$.
\item $\int_{\R}aq_0(\tilde{w})q_1(w(0),\tilde{w})\diff x+\int_{\R}a'\tilde{q}_0(\tilde{w})q_2(w(0),\tilde{w})\diff x=0$.
\item For every $1 \leq i \leq 6$, the gradient of $q_i$ is bounded on $V \times \tilde{V}$ and the gradients of $q_0$ and $\tilde{q}_0$ are bounded on $\tilde{V}$.
\end{enumerate}
Denote $E_{K,T^*}:=\lbrace X \in C^1([0,T^*]):X(0)=0 \,$and $\|\dot{X}\|_{\infty} \leq K\rbrace$.
Then, for any $\delta, K > 0$, there exists $T^*(K, C, \delta)>0$ such that:
\begin{align*}
\sup_{X \in E_{K,T^*}}\sup_{0 \leq t \leq T^*}\Bigl|&\dot{X}(t)\left(\int_{\R}a(q_0(\tilde{w}))_x\cdot q_1(w^{X(t)}(t),\tilde{w})\diff x+\int_{\R}a'q_2(w^{X(t)}(t),\tilde{w})\diff x\right) \\ &+
\int_{\R}a\tilde{q}_0(\tilde{w})q_3(w^{X(t)}(t),\tilde{w})\diff x+\int_{\R}a''q_4(w^{X(t)}(t),\tilde{w})\diff x \\&-
\int_{\R}a\left(q_5(w^{X(t)}(t),\tilde{w})\right)_x^2\diff x+\int_{\R}a'q_6(w^{X(t)}(t),\tilde{w})\diff x \\
&- \dot{X}(0)\left(\int_{\R}a(q_0(\tilde{w}))_xq_1(w(0),\tilde{w})\diff x+\int_{\R}a'q_2(w(0),\tilde{w})\diff x\right) \\
&- \int_{\R}a\tilde{q}_0(\tilde{w})q_3(w(0),\tilde{w})\diff x-\int_{\R}a''q_4(w(0),\tilde{w})\diff x\\
&+ \int_{\R}a\left(q_5(w(0),\tilde{w})\right)_x^2\diff x-\int_{\R}a'q_6(w(0),\tilde{w})\diff x\Bigr| \leq \delta.
\end{align*}
\end{lemma}
We apply with the lemma with the following functions (cf. \eqref{step3}):
\begin{enumerate}
    \item $q_0(v)=v$,
    \item $\tilde{q}_0(v)=v$,
    \item $q_1(u,v)=-2(u-v)$,
    \item $q_2(u,v)=(u-v)^2$,
    \item $q_3(u,v)=(-2\sigma(u-v)-2A(u|v))$,
    \item $q_4(u,v)=(u-v)^2$,
    \item $q_5(u,v)=\sqrt{2}(u-v)$,
    \item $q_6(u,v)=A'(v)(u-v)^2+\frac{1}{2}\left(\int_v^uG^{(3)}(k)-2SA^{(3)}(k)(u-k)^2dk\right)$.
\end{enumerate}
So, there exists $T^*$ positive such that for all $0 \leq s \leq T^*$ and $n$ in $\N$
\[
\dv{}{t}\int_\mathbb{R} a(x)|u(x+X_n(t),t)-S(x)|^2\diff x|_{t=s} \geq \frac{c}{2}.
\]
Then, integrating from $0$ to $T^*$, we find that for any $0 < t \leq T^*$ and $n \in \mathbb{N}$:
\begin{align*}
    \int_\mathbb{R} a(x)|u(x+X_n(t),t)-S(x)|^2\diff x \geq \dfrac{ct}{2}+\int_\mathbb{R} a(x)|u_0(x)-S(x)|^2\diff x 
\end{align*}
We may take the limit in $n$ (utilizing the fact that $\|X_n-X\|_\infty \to 0$) to obtain for all $t \in (0,T^*]$
\[
\int_\mathbb{R} a(x)|u(x+X(t),t)-S(x)|^2\diff x > \int_\mathbb{R} a(x)|u_0-S(x)|^2\diff x,
\]
as desired.

\section{Preliminaries for the barotropic Navier-Stokes case}
For the latter half of the paper, we will show the failure of $a$-contraction for sufficiently large shocks of the barotropic Navier-Stokes system. The structure will be the same as in the scalar case: we will discuss existence/uniqueness of perturbations, the relative entropy structure, the linearization of the shift condition, and the construction of the initial perturbation, and finally the approximation of the inital perturbation and shifts, in the same order. \\ 
Let $\gamma>1$. Fix $\tilde{a}$ in $W^{2,\infty}([0,1])$ such that $\tilde{a}$ is valued in $\R^+$.
\subsection{Transformation of system \eqref{bns}}
We make a change of variables in \eqref{bns} using the effective velocity, $h:=u+p(v)_x$ (this was first considered by Shelukhin \cite{MR0832100} and generalized by Desjardins \cite{MR1989675}, \cite{MR2257849}, \cite{MR1978317} and Haspot \cite{MR3013413}, \cite{MR3620595}. See also \cite{MR4195742}). Under this transformation, the system \eqref{bns} becomes:
\begin{equation}\label{transformedbns}
\begin{cases}
v_t-h_x=-p(v)_{xx}, & (x,t) \in \mathbb{R} \times (0,\infty), \\ 
h_t+p(v)_x=0, & (x,t) \in \mathbb{R} \times (0,\infty), \\
(v(x,0),h(x,0))=(v_0(x), h_0(x)), & x \in \R.
\end{cases}
\end{equation}
The transformed viscous shock $(\tilde{v}, \tilde{h})$ satisfies
\begin{align}\label{transformedviscousshock}
\begin{cases}
-\sigma \tilde{v}'-\tilde{h}'=-p(\tilde{v})'', \\
-\sigma \tilde{h}'+p(\tilde{v})'=0, \\
\lim_{\xi \to \pm \infty}=(\tilde{v}, \tilde{h})=(v_\pm, u_\pm), \\
\lim_{\xi \to \pm \infty}(\tilde{v}', \tilde{h}')=(0,0), \\
\sigma=-\sqrt{-\frac{p(v_+)-p(v_-)}{v_+-v_-}}, \text{ if } v_->v_+,\\
\sigma=\sqrt{-\frac{p(v_+)-p(v_-)}{v_+-v_-}}, \text{ if } v_-<v_+, \\
\tilde{v}'=-\frac{-\sigma(\tilde{v}-v_-)-\frac{p(\tilde{v})-p(v_-)+\sigma u_-}{\sigma}+u_-}{p'(\tilde{v})}=:q(\tilde{v}). \\
\end{cases}
\end{align}
where $\tilde{h}=\tilde{u}+p(\tilde{v})_x$. The viscous shock exists if and only if the endstates $(v_-,u_-)$ and $(v_+,u_+)$ satisfy the conditions \eqref{rh}. The existence of the wave is proven in \cite{MR2785875} by reducing the problem to a scalar differential equation. \\
Notice that for any $(v_-,u_-)$ in $\R^+ \times \R$, we have that for every $v_+$ in $(0,v_-)$, there exists a unique $u_+$ in $\R$ and a unique $\sigma$ in $\R^-$ such that $(\sigma,v_-,u_-,v_+,u_+)$ satisfies the conditions \eqref{rh}. In the following, we fix $(v_-,u_-)$ in $\R^+ \times \R$ and consider $v_+$ in $(0,v_-)$ as a variable that will be chosen close enough to $0$. \\
For convenience, we work in the moving frame by making the change of variable $(x,t) \mapsto (\xi=x-\sigma t,t)$. Then, \eqref{transformedbns} becomes:
\begin{equation}\label{transformedbnsmf}
\begin{cases}
v_t-\sigma v_\xi-h_\xi=-p(v)_{\xi \xi}, & (\xi,t) \in \mathbb{R} \times (0,\infty), \\ 
h_t-\sigma h_\xi+p(v)_\xi=0, & (\xi,t) \in \mathbb{R} \times (0,\infty), \\
(v(x,0),h(x,0))=(v_0(x), u_0(x)), & \xi \in \R.
\end{cases}
\end{equation}
Theorem \ref{maintheorem2} follows directly from the following analogous theorem for the transformed system \eqref{transformedbnsmf}, which shows a failure of $a$-contraction with respect to the entropy $\eta=\frac{h^2}{2}+Q(v)$, where $Q'(v)=-p(v)$:
\begin{theorem}\label{maintheorem3}
Consider the system \eqref{transformedbnsmf}. Fix $\tilde{a} \in W^{2,\infty}([0,1])$. Let $(v_-,u_-) \in \mathbb{R}^+ \times \mathbb{R}$ be a given constant state. There exists some $\epsilon$ positive such that the following is true. Let $(v_+,u_+) \in \mathbb{R}^+ \times \mathbb{R}$ satisfy \eqref{rh} and $|v_+| < \epsilon$. Then, for every $\delta > 0$ and $s > 0$, there exists smooth initial data $(v_0,u_0)$ to \eqref{transformedbnsmf} with $\tilde{U}-(u_0,u_0) \in (C_0^\infty(\mathbb{R}))^2$ (where $\tilde{U}=(\tilde{v},\tilde{h})$ denotes the viscous shock connecting $(v_-,u_-), (v_+,u_+)$), such that for any Lipschitz shift $X$ and every $t \in (0, T^*)$:
\[
\int_\mathbb{R}a(\xi)\eta(u(\xi+X(t),t)|\tilde{U})\diff \xi > \int_\mathbb{R}a(\xi)\eta((v_0(\xi),u_0(\xi)|\tilde{U})\diff \xi,
\]
for some $T^*=T^*(\tilde{a},\|\dot{X}\|_\infty)$, where $U$ is the solution to \eqref{transformedbnsmf}.
\end{theorem}

\subsection{Existence and uniqueness of perturbations}
For the existence and uniqueness of solutions, we refer to the article \cite{MR4113341} in which smooth solutions to the system in Eulerian coordinates are given for initial perturbations in some Sobolev spaces (which is equivalent to the Lagrangian mass formulation for smooth enough solutions). As our initial perturbations are in $C_0^\infty$, this result provides existence and uniqueness of global-in-time smooth solutions. Again, bounds on solutions which are initially perturbations of a viscous shock wave will be discussed before the approximation of shifts in Section \ref{systemshiftapprox}.

\subsection{Relative entropy structure}
To calculate the time derivative of the quantity $\int_\mathbb{R}a(x)\eta(u^{-X}|\tilde{u})\diff x$, we refer to the article \cite{MR4195742}, in which the same calculation was made. It is identical (but in a different setting) to the calculation made in Lemma \ref{structurelemma}. Recall that in this case, we have the entropy function:
\begin{align*}
\eta(v,h)&=\frac{h^2}{2}+Q(v), \notag \\
Q'(v)&=-p(v), \notag \\
\end{align*}
and the relative quantities $\eta(u|\tilde{u}), Q({v|\tilde{v}}), p(v|\tilde{v})$ are defined analogously as in Section \ref{relentropyscalar}. Recall that we have defined:
\begin{align*}
    b(\tilde{v}):=\frac{p(\tilde{v})-p(v_-)}{p(v_+)-p(v_-)}.
\end{align*}
\begin{remark} As in the scalar case, we find it convenient to work in the shock variables. For this problem, the change of variables we choose is $y=\tilde{v}(\xi)$, and, given the shift function $X:\mathbb{R}_{\geq 0} \to \mathbb{R}$, we define the quantities:
\begin{align*}
    &w(y,t)=v^{-X}(\tilde{v}^{-1}(y),t)-y, \\
    &g(y,t)=h^{-X}(\tilde{v}^{-1}(y),t)-\tilde{h}(\tilde{v}^{-1}(y)).
\end{align*}
For convenience, we also have the following notation:
\[
[p]:=p(v_+)-p(v_-).
\]
Note that for $v_-$ fixed, $[p] \to \infty$ as $v_+ \to 0$.
\end{remark}
Then, we obtain the following:
\begin{lemma}\label{systemstructurelemma}
Let $\tilde{U}=(\tilde{v},\tilde{h})$ be the viscous shock solution defined by \eqref{transformedviscousshock}. Let $U=(v,h)$ be the solution to \eqref{transformedbnsmf}. Let $a(\xi)=\tilde{a}(b((\tilde{v}(\xi)))$ for some $\tilde{a} \in W^{2,\infty}([0,1])$. 
Then, we have:
\begin{align}\label{systemrepresentation}
\dv{}{t}\int_\mathbb{R}a(x)\eta(U^{-X}|\tilde{U})\diff \xi=-\dot{X}(t)\tilde{Y}((w,g))+\tilde{B}(w)-\tilde{G}((w,g)),
\end{align}
where
\begin{align*}
    \tilde{Y}((w,g))=&\int_{v_+}^{v_-}\tilde{a}'(b(y))\frac{p'(y)}{[p]}Q(w+y|y)\diff y+\int_{v_+}^{v_-}\tilde{a}'(b(y))\frac{p'(y)}{[p]}\frac{g(y)^2}{2}\diff y \\
    &+\int_{v_+}^{v_-}\tilde{a}(b(y))p'(y)w(y)\diff y- \frac{1}{\sigma}\int_{v_+}^{v_-}\tilde{a}(b(y))p'(y)g(y)\diff y. \notag
\end{align*}
Further:
\begin{align*}
    \tilde{B}(w)=& -\frac{1}{2[p]\sigma}\int_{v_+}^{v_-}\tilde{a}'(b(y))p'(y)|p(w+y)-p(y)|^2\diff y \notag \\
    &-\sigma \int_{v_+}^{v_-}\tilde{a}(b(y))p(w+y|y)\diff y \notag \\
    &-\frac{1}{2[p]^2}\int_{v_+}^{v_-}\tilde{a}''(b(y))(p'(y))^2q(y)|p(w+y)-p(y)|^2\diff y \notag \\
    &-\frac{1}{2[p]}\int_{v_+}^{v_-}\tilde{a}'(b(y))(\sigma+\frac{p'(y)}{\sigma})|p(w+y)-p(y)|^2\diff y, \notag \\
\end{align*}
and:
\begin{align*}
    \tilde{G}((w,g))=&-\frac{\sigma}{2[p]}\int_{v_+}^{v_-}\tilde{a}'(b(y))p'(y)\left(g(y)-\frac{p(w+y)-p(y)}{\sigma}\right)^2\diff y \notag \\
    &-\frac{\sigma}{[p]} \int_{v_+}^{v_-}\tilde{a}'(b(y))p'(y)Q(w+y|y)\diff y \notag \\
    &-\int_{v_+}^{v_-}\tilde{a}(b(y))|p'(w+y)w_y(y)+p'(w+y)-p'(y)|^2q(y)\diff y.
\end{align*}
\end{lemma}
\begin{proof}
By Lemma 2.3 in \cite{MR4195742}, we immediately have:
\begin{align}\label{systemrepresentationkv}
\dv{}{t}\int_\mathbb{R}a(x)\eta(U^{-X}|\tilde{U})\diff \xi=\dot{X}(t)Y(U^{-X})+B(U^{-X})-G(U^{-X}),
\end{align}
where:
\begin{align*}
    Y(U):=&\int_\mathbb{R}a'\eta(U|\tilde{U})\diff \xi-\int_\mathbb{R}a(\nabla \eta(\tilde{U}))_\xi \cdot (U-\tilde{U})\diff \xi, \notag \\
    B(U):=&\frac{1}{2\sigma}\int_\mathbb{R}a'|p(v)-p(\tilde{v})|^2\diff \xi+\sigma \int_\mathbb{R}a\tilde{v}_\xi p(v|\tilde{v})\diff \xi +\frac{1}{2}\int_\mathbb{R}a''|p(v)-p(\tilde{v})|^2\diff \xi, \notag \\
    G(U):=&\frac{\sigma}{2}\int_\mathbb{R}a'\left(h-\tilde{h}-\frac{p(v)-p(\tilde{v})}{\sigma}\right)^2\diff \xi \\
    &+\sigma\int_\mathbb{R}a'Q(v|\tilde{v})\diff \xi+\int_\mathbb{R}a|(p(v)-p(\tilde{v}))_\xi|^2\diff \xi.
\end{align*}
Making a change of variables $y=\tilde{v}(\xi)$ in \eqref{systemrepresentationkv}, we obtain the result. 
\end{proof}

\subsection{Reduction of the shift condition}
As before, we need to construct an initial perturbation $(w,g)$ with $\tilde{Y}((w,g))=0$. We again show that it suffices to construct a perturbation that fits a certain linear condition.
\begin{lemma}\label{systemift}
Let $(w,g) \in (C_0^\infty([v_+,v_-]))^2$ such that \newline 
$\int_{v_+}^{v_-}\tilde{a}(b(y))p'(y)w(y)-\frac{1}{\sigma}\int_{v_+}^{v_-}\tilde{a}(b(y))p'(y)g(y)\diff y=0$. Let $\phi \in C_0^\infty([v_+,v_-])$ such that $\int_{v_+}^{v_-}\tilde{a}(b(y))p'(y)\phi(y) \diff y \neq 0$. Denote 
\[
\lambda_*:=\frac{\int_{v_+}^{v_-}\tilde{a}'(b(y))\frac{p'(y)^2}{[p]}w(y)^2\diff y-\int_{v_+}^{v_-}\tilde{a}'(b(y))\frac{p'(y)}{2[p]}g(y)^2\diff y}{\int_{v_+}^{v_-}\tilde{a}(b(y))p'(y)\phi(y) \diff y} \, .
\]
There exists some $\epsilon_0$ positive, some $\delta$ positive and some $\lambda : (0,\epsilon_0) \rightarrow (\lambda_*-\delta,\lambda_*+\delta)$, such that for every $\epsilon$ in $(0,\epsilon_0)$, with $(\overline{w}, \overline{g})=\epsilon(w,g)+(\epsilon^2\lambda(\epsilon) \phi,0)$ we have that
\[
\tilde{Y}((\overline{w}, \overline{g}))=0 \,,
\]
and such that
\begin{align*}
    \lim_{ \epsilon \rightarrow 0 }\lambda(\eps)=\lambda_* \, .
\end{align*}
\end{lemma}
The proof is similar to that of Lemma \ref{ift}, so we focus on the new features. Notice that here, it suffices to make a small perturbation in one variable only. 
\begin{proof}
    Let $\lambda_*$ be as above. \\
    With $(\overline{w},\overline{g})$ as above and $\eps_1$ and $\delta_1$ small enough such that there exists $c$ positive such that for every $\eps$ in $(-\eps_1,\eps_1)$, $\lambda$ in $(\lambda_*-\delta_1,\lambda_*+\delta_1)$ and $x$ in $\R$, $\tilde{v}(x)+\eps w(\tilde{v}(x))+\eps \lambda \phi(\tilde{v}(x)) \geq c$. Then, for every $\eps$ in $(-\eps_1,\eps_1)$ and $\lambda$ in $(\lambda_*-\delta_1,\lambda_*+\delta_1)$
    \begin{align*}
        \tilde{Y}(\overline{w},\overline{g})=&\int_{v_+}^{v_-}\tilde{a}'(b(y))\frac{p'(y)}{[p]}Q(\overline{w}+y|y)\diff y+\eps^2\int_{v_+}^{v_-}\tilde{a}'(b(y))\frac{p'(y)}{[p]}\frac{g^2}{2}\diff y \\
        &+\int_{v_+}^{v_-}\eps^2\tilde{a}(b(y))p'(y)\phi(y)\lambda \diff y \, ,
    \end{align*}
    and notice that for every $y$ in $(v_+,v_-)$
    \begin{align*}
        Q(w+y|y)=-\eps^2\int_0^1p'(y+\theta \overline{w}(y))(1-\theta)\diff\theta (w(y)+\eps \lambda \phi(y))^2 \, .
    \end{align*}
    Now, consider
    \begin{align*}
        f(\eps,\lambda)=&-\int_{v_+}^{v_-}\int_0^1\tilde{a}'(b(y))\frac{p'(y)}{[p]}p'(y+\eps\theta (w(y)+\eps \lambda \phi(y)))(1-\theta)\diff\theta (w(y) + \eps \lambda \phi(y))^2\diff y \\
        &+\int_{v_+}^{v_-}\tilde{a}'(b(y))\frac{p'(y)}{[p]}\frac{g^2}{2}\diff y+\lambda\int_{v_+}^{v_-}\tilde{a}(b(y))p'(y)\phi(y) \diff y \, ,
    \end{align*}
    defined for $(\eps,\lambda)$ in $(-\eps_1,\eps_1)\times(\lambda_*-\delta,\lambda_*+\delta)$. \\
    We then notice that $f(0,\lambda_*)=0$ and that by smoothness of $p$ and derivation under the integral sign, $f$ is smooth in $(\eps,\lambda)$ and that
    \begin{align*}
        D_{\lambda}f(0,\lambda_*)=\int_{v_+}^{v_-}\tilde{a}'(b(y))p'(y)\phi(y)\diff y \, .
    \end{align*}
    Hence, we can now about the existence of such a $\lambda:[0,\eps_0] \rightarrow (\lambda_*-\delta,\lambda_*+\delta)$ such that for every $\eps$ in $[0,\eps_0]$, $f(\eps,\lambda(\eps))=0$.
\end{proof}
As before, note that we may choose $\epsilon$ sufficiently small such that the $\lambda$ given by Lemma \ref{systemift} is arbitrarily close to $\lambda_*$.
\subsection{Estimates for various quantities}
In this section, we provide estimates for various quantities depending on $v_-$ and $v_+$ that will be needed later. These will be used crucially to determine the sign of the quadratic form considered later on.
\begin{lemma}\label{alphaestimate}
There exists positive constants $C$ and $\delta$ such that for any $v_+$ in $\R^+$ such that $v_+<\min(\delta,v_-)$ the following holds. \\ 
Define $\alpha=\frac{\int_{v_+}^{v_-}\tilde{a}(b(y))\sigma \diff y}{\int_{v_+}^{2v_+}\tilde{a}(b(y))p'(y)\diff y}$, with $\sigma=-\sqrt{\frac{p(v_+)-p(v_-)}{v_--v_+}}$. Then, we have:
\begin{equation}
|\alpha| \leq C\frac{1}{\sqrt{p(v_+)}}.
\end{equation}
\end{lemma}
\begin{proof}
Using the fact that $\tilde{a}$ is continuous on $[0,1]$ and only takes positive values, we have that for every $v_-$ and $v_+$ in $\R^+$ such that $v_+<v_-$
\begin{equation*}
|\alpha| \leq C \frac{\sigma(v_--v_+)}{p(2v_+)-p(v_+)}.
\end{equation*}
Now, recalling \eqref{transformedviscousshock}, we obtain:
\begin{equation*}
|\alpha| \leq C\frac{\sqrt{v_-}\sqrt{p(v_+)}}{(1-2^{-\gamma})p(v_+)}= C\frac{\sqrt{v_-}}{(1-2^{-\gamma})\sqrt{p(v_+)}}.
\end{equation*}
\end{proof}
Recall from \eqref{transformedviscousshock} that $q$ is defined
\begin{equation*}
q(y)=-\frac{-\sigma(y-v_-)-\frac{p(y)-p(v_-)+\sigma u_-}{\sigma}+u_-}{p'(y)}.
\end{equation*}
\begin{lemma}\label{qestimate}
There exist positive constants $C$ and $\delta$ such that for every $v_+$ in $\R^+$ such that $v_+<\min(\delta,v_-)$, we have that
\begin{equation}
|q(y)| \leq C\frac{|\sigma|}{|p'(y)|} \, .
\end{equation}
\end{lemma}
\begin{proof}
We compute:
\begin{equation*}
|q(y)| \leq \frac{|\sigma|v_-+C\frac{p(v_+)}{|\sigma|}}{|p'(y)|} \, .
\end{equation*}
Recalling from $\eqref{transformedviscousshock}$ that $|\sigma| \leq C\sqrt{p(v_+)}$ as $v_+ \to 0$ (for $v_-$ fixed), we get:
\begin{equation*}
|q(y)| \leq C\frac{|\sigma|}{|p'(y)|} \leq C\frac{|\sigma|}{|p'(y)|} \, .
\end{equation*}
\end{proof}
\section{Proof of Theorem \ref{maintheorem3}}
In this section, we prove Theorem \ref{maintheorem3}. As in the scalar case, the strategy is to make a small initial perturbation that causes the derivative of the relative entropy given by the representation \eqref{systemrepresentation} to be positive at $t=0$.

\begingroup
\allowdisplaybreaks
\subsection{Small perturbation of $\tilde{U}$}
Firstly, assume that we can construct a perturbation $(w,g)$ with $\int_{v_+}^{v_-}\tilde{a}(b(y))p'(y)w(y)-\frac{1}{\sigma}\int_{v_+}^{v_-}\tilde{a}(b(y))p'(y)g(y)\diff y=0$. Then, for $\epsilon < \epsilon_0$, consider the perturbation associated with $(\overline{w}, \overline{g})$ constructed in Lemma \ref{systemift} (that is, $u_0(x)=(\tilde{v}(x)+\overline{w}(\tilde{v}(x)),\tilde{h}(x)+\overline{g}(\tilde{v}(x)))$). Plugging this into \eqref{systemrepresentation}, and using a Taylor expansion of all relative quantities (and also quantities of the form $p(w+y)-p(y)$), using the integral form of Taylor remainder to control the higher order terms, we obtain, as in Section \ref{smallperturbation}, for fixed $v_-, v_+$:
\begin{align}\label{eperturbationsystem}
\dv{}{t}\int_\mathbb{R}a(\xi)\eta(u^{-X}|\tilde{u})\diff \xi|_{t=0}=\epsilon^2 \Bigg(&-\frac{1}{2[p]\sigma}\int_{v_+}^{v_-}\tilde{a}'(b(y))(p'(y))^3w(y)^2\diff y  \notag \\
    &-\frac{\sigma}{2} \int_{v_+}^{v_-}\tilde{a}(b(y))p''(y)w(y)^2\diff y \notag\\ 
    &-\frac{1}{2[p]^2}\int_{v_+}^{v_-}\tilde{a}''(b(y))(p'(y))^4q(y)w(y)^2\diff y \notag \\ &-\frac{1}{2[p]}\int_{v_+}^{v_-}\tilde{a}'(b(y))\left(\sigma+\frac{p'(y)}{\sigma}\right)(p'(y))^2w(y)^2\diff y \notag \\
    &+\frac{\sigma}{2[p]}\int_{v_+}^{v_-}\tilde{a}'(b(y))p'(y)\left(g(y)-\frac{p'(y)w(y)}{\sigma}\right)^2\diff y \notag \\
    &-\frac{\sigma}{2[p]} \int_{v_+}^{v_-}\tilde{a}'(b(y))(p'(y))^2w(y)^2\diff y \notag \\
    &+\int_{v_+}^{v_-}\tilde{a}(b(y))((p'w)_y(y))^2q(y)\diff y+R(v_+,\eps)\Bigg),
\end{align} 
where:
\begin{align*}
R(v_+,\epsilon)=&\frac{1}{\eps^2}\left(\tilde{B}(\overline{w})-\tilde{G}((\overline{w}, \overline{g}))\right)-\frac{1}{2[p]\sigma}\int_{v_+}^{v_-}\tilde{a}'(b(y))(p'(y))^3w(y)^2\diff y  \notag \\
&-\frac{\sigma}{2} \int_{v_+}^{v_-}\tilde{a}(b(y))p''(y)w(y)^2\diff y \notag\\ 
&-\frac{1}{2[p]^2}\int_{v_+}^{v_-}\tilde{a}''(b(y))(p'(y))^4q(y)w(y)^2\diff y \notag \\ &-\frac{1}{2[p]}\int_{v_+}^{v_-}\tilde{a}'(b(y))\left(\sigma+\frac{p'(y)}{\sigma}\right)(p'(y))^2w(y)^2\diff y \notag \\
&+\frac{\sigma}{2[p]}\int_{v_+}^{v_-}\tilde{a}'(b(y))p'(y)\left(g(y)-\frac{p'(y)w(y)}{\sigma}\right)^2\diff y \notag \\
&-\frac{\sigma}{2[p]} \int_{v_+}^{v_-}\tilde{a}'(b(y))(p'(y))^2w(y)^2\diff y \notag \\
&+\int_{v_+}^{v_-}\tilde{a}(b(y))(p'w)_y(y)^2q(y)\diff y \, .
\end{align*}
\endgroup
For every $v_+$ in $(0,v_-)$ there exists positive constants $C$ and $\eps_0$ such that for every $\eps$ in $(0,\eps_0)$, we have that $|R(v_+,\eps)| \leq C\eps$ (by using Taylor expansion with integral remainder).

\subsubsection{Construction of $w$, $g$}
Define the following functional $F$:
\begin{equation*}
    F((w,g))):=\sum_{i=1}^3J_i((w,g)) \, ,
\end{equation*}
where:
\begin{align*}
    J_1((w,g)):=& \frac{-\sigma}{2[p]} \int_{v_+}^{v_-}\tilde{a}'(b(y))(p'(y))^2w(y)^2\diff y \notag \\
    &+\int_{v_+}^{v_-}\tilde{a}(b(y))((p'w)_y)(y)^2q(y)\diff y \, , \notag \\
    J_2((w,g)):=&-\frac{\sigma}{2} \int_{v_+}^{v_-}\tilde{a}(b(y))p''(y)w(y)^2\diff y \, , \notag \\
    J_3((w,g)):=&-\frac{1}{2[p]\sigma}\int_{v_+}^{v_-}\tilde{a}'(b(y))(p'(y))^3w(y)^2\diff y \notag \\
    &-\frac{1}{2[p]^2}\int_{v_+}^{v_-}\tilde{a}''(b(y))(p'(y))^4q(y)w(y)^2\diff y \notag \\
    &-\frac{1}{2[p]}\int_{v_+}^{v_-}\tilde{a}'(b(y))(\sigma+\frac{p'(y)}{\sigma})(p'(y))^2w(y)^2\diff y \notag \\
    &+\frac{\sigma}{2[p]}\int_{v_+}^{v_-}\tilde{a}'(b(y))p'(y)\left(g(y)-\frac{p'(y)w(y)}{\sigma}\right)^2\diff y \notag \, . \\
\end{align*}
Due to the representation \eqref{eperturbationsystem}, it suffices to show that for $v_+$ sufficiently close to zero, we may find a pair of smooth functions on $[v_+,v_-]$ such that $F((w,g)) > 0$. Consider the following functions:
\begin{align*}
    w(y)&=\frac{1}{p'(y)}, \\
    g(y)&=\alpha \mathds{1}_{[v_+,2v_+]}(y) \, .
\end{align*}
(recall the definition of $\alpha$ from Lemma \ref{alphaestimate}). Note that we have picked $\alpha$ exactly so that $(w,g)$ satisfy the hypotheses of Lemma \ref{systemift}. We will show that this pair $(w,g)$ satisfies $F(w,g) > 0$ for $v_+$ sufficiently small. Then, we will approximate $w$ and $g$ by smooth functions that are zero on the boundary in the topology induced by $F$, as in the scalar case. \\
To do so, we will control the different terms here by considering a positive constant $C$, possibly changing from line to line. Crucially, it will not depend on $v_+$ but only on the other variables of the problem.
\subsubsection{Control on $J_1$}
There exists positive constants $C$ and $\delta$ such that for every $v_+$ in $(0,\min(\delta,v_-))$
\begin{equation}\label{J_1estimate}
\left|J_1((w,g))\right|=\left|\frac{-\sigma}{2[p]} \int_{v_+}^{v_-}\tilde{a}'(b(y))\diff y\right| \leq C\frac{|\sigma|}{[p]} \, .
\end{equation}
as $(p'w)_y(y)=0$.

\subsubsection{Control on $J_2$}
This is the positive term. There exists positive constants $c,c'$ and $\delta$ such that for every $v_+$ in $(0,\min(\delta,v_-))$
\begin{equation}\label{J_2estimate}
J_2((w,g))=-\sigma \int_{v_+}^{v_-}\tilde{a}(b(y))p''(y)w(y)^2\diff y \geq -\sigma c'\int_{v_+}^{v_-}y^\gamma \diff y \geq -\sigma c \, ,
\end{equation}
(recall that $\sigma < 0$, as we consider $1$-shocks), and $c$ is a positive constant independent of $v_+$.

\subsubsection{Control on $J_3$}
We handle the four terms in $J_3$ separately. There exists positive constants $C$ and $\delta$ such that for every $v_+$ in $(0,\min(\delta,v_-))$
\begin{equation}\label{J_3_1estimate}
\left|-\frac{1}{2[p]\sigma}\int_{v_+}^{v_-}\tilde{a}'(b(y))(p'(y))^3w(y)^2\diff y\right| \leq \frac{C}{[p]\left|\sigma\right|}\int_{v_+}^{v_-}p'(y)\diff y \leq \frac{C}{\left|\sigma\right|} \leq C\frac{\left|\sigma\right|}{[p]} \, .
\end{equation}
For the second, recalling Lemma \ref{qestimate} for the bound on $q$:
\begin{equation}\label{J_3_2estimate}
\left|\frac{1}{2[p]^2}\int_{v_+}^{v_-}\tilde{a}''(b(y))(p'(y))^4q(y)w(y)^2\diff y\right| \leq \frac{C\left|\sigma\right|}{[p]} \, .
\end{equation}
For the third, we check:
\begin{equation}\label{J_3_3estimate}
\left| \frac{1}{2[p]}\int_{v_+}^{v_-}\tilde{a}'(b(y))(\sigma+\frac{p'(y)}{\sigma})(p'(y))^2w(y)^2\diff y \right| \leq \frac{C|\sigma|}{[p]} \, .
\end{equation}
Finally, for the last, we break the integral into two parts and calculate each separately:
\begin{align*}
 &\left|\frac{\sigma}{2[p]}\int_{v_+}^{v_-}\tilde{a}'(b(y))p'(y)\left(g(y)-\frac{p'(y)w(y)}{\sigma}\right)^2\diff y\right| \\
 &\leq \left|\frac{\sigma}{2[p]}\int_{v_+}^{2v_+}\tilde{a}'(b(y))p'(y)\left(\alpha-\frac{1}{\sigma}\right)^2\diff y\right| \notag \\
 & \quad +\left|\frac{\sigma}{2[p]}\int_{2_+}^{v_-}\tilde{a}'(b(y))p'(y)\left(\frac{1}{\sigma}\right)^2\diff y\right| \, .
\end{align*}
For the first term, recalling Lemma \ref{alphaestimate} and \eqref{transformedviscousshock}, we see that $(\alpha-\frac{1}{\sigma})$ grows asymptotically like $\frac{1}{\sqrt{p(v_+)}}$. So, we obtain:
\begin{equation*}
\left|\frac{-\sigma}{2[p]}\int_{v_+}^{2v_+}\tilde{a}'(b(y))p'(y)\left(\alpha-\frac{1}{\sigma}\right)^2\diff y\right| \leq \left|\frac{C\sigma}{2[p]}\frac{1}{p(v_+)}(p(2v_+)-p(v_-))\right| \leq \left|\frac{C\sigma}{[p]}\right| \, .
\end{equation*}
For the second, we get:
\begin{equation*}
\left|\frac{-\sigma}{2[p]}\int_{2_+}^{v_-}\tilde{a}'(b(y))p'(y)\left(\frac{1}{\sigma}\right)^2\diff y\right| \leq \left|\frac{C}{\sigma}\right| \, .
\end{equation*}

Finally, consolidating all the previous estimates, we get that there exists positive constant $\delta$ and $C$ such that for every $v_+$ in $(0,\min(\delta,v_-))$
\begin{equation}\label{J_3estimate}
\left|J_3((w,g))\right| \leq \left|\frac{C\sigma}{[p]}\right| \, .
\end{equation}
\subsubsection{Conclusion}
Hence, there exists $\delta$ positive such that for every $v_+$ in $(0,\min(\delta,v_-))$, we have that
\begin{align*}
    F((w,g)) \geq \frac{|\sigma|c}{2} \, ,
\end{align*}
with $w$ and $g$ as above. Further, $(w,g)$ satisfies the linear shift condition in Lemma \ref{systemift}.
\subsection{Approximation in the functional topology}
Now, we again show that we may find functions $(\tilde{w}, \tilde{g}) \in [C_0^\infty([v_+,v_-])]^2$ with $F((\tilde{w},\tilde{g})) > 0$, such that the linear shift condition is satisfied. For $v_+, v_-$ fixed, define the following norm:
\[
\|(b,m)\|_H^2:=\|b\|_2^2+\|m\|_2^2+\int_{v_+}^{v_-}(y-v_+)(y-v_-)|m_y(y)|^2\diff y.
\]
This is the $L^2$ norm in $b$, and the same norm used in Section \ref{initdataconstruct} for the norm in $m$. Define $H$ to be the closure of $(C_0^\infty((v_+,v_-))^2$ in this topology. Note that $F$ is continuous on $H$ equipped with this topology. Further, we have that $(w,g) \in H$; indeed, doing the same approximation as in Section \ref{initdataconstruct} suffices. Finally, by invoking Lemma \ref{kernellemma} again, we may construct $(\tilde{w},\tilde{g}) \in (C_0^\infty((v_+,v_-)))^2$ such that the linear shift condition is satisfied, and $F((\tilde{w},\tilde{g})) > 0$, as desired.

\subsection{Approximation of Lipschitz shifts}\label{systemshiftapprox}
Now, we show the approximation of a Lipschitz shift by a class of $C^1$ shifts. The strategy is identical is done in the scalar case. For the approximating shifts, we will use the same approximation as in the scalar case, invoking Lemma \ref{c1shifts}. For the regularity on $U-\tilde{U}$, we consider the following lemma:
\begin{lemma}
Let $\tilde{U}=(\tilde{v},\tilde{h})$ be the viscous shock solution defined by \eqref{transformedviscousshock}. Let $v_0$ and $h_0$ be such that $(v_0-\tilde{v}, h_0-\tilde{h}) \in C_0^\infty(\mathbb{R})^2$. Let $U=(v,h)$ be the classical solution to \eqref{transformedbnsmf} associated with the initial data $(v_0,h_0)$. Then, there exists some $T$ positive such that $(v-\tilde{v}, h-\tilde{h})$ is in $C^0([0,T]; H^2(\mathbb{R})^2)$. Furthermore, there are uniform upper and lower bounds on $v$ valid for every $t$ in $[0,T]$ and every $\xi$ in $\mathbb{R}$ of the form $\underline{\kappa}(T) \leq v(\xi,t) \leq \overline{\kappa}(T) \ $. 
\end{lemma}
We will only sketch the proof of this lemma.
\begin{proof}
For initial perturbations in $H^k$ (for $k \geq 3$), we can prove the existence and uniqueness of classical solutions by some usual scheme. \\
The existence relies on a viscous approximation. \\
More details for perturbations of a constant state but in more generality can be found in \cite{MR2731999}. Here, we consider perturbation of a nonconstant state, hence extra terms arise, but the method applies.\\
Denote 
\begin{align*}
    A:&\R \times \R \times \R \rightarrow \mathcal{M}_2(\R) \\
    &(v,h) \mapsto \begin{pmatrix}
    \sigma-p''(v)\tilde{v}'(x) & 1 \\
    -p'(v) & \sigma
    \end{pmatrix} \, , \\
    B:&\R^+ \times \R \rightarrow \mathcal{M}_2(\R) \\
    &(v,h) \mapsto \begin{pmatrix}
    -p'(v) & 0 \\
    0 & 0
    \end{pmatrix} \, , \\
    C:&\R^+ \times \R \times \R \rightarrow \mathcal{M}_2(\R) \\
    &(v,h,x) \mapsto \begin{pmatrix}
-\int_0^1p''(\tilde{v}(x)+\theta v)\tilde{v}''(x)+p^{(3)}(\tilde{v}(x)+\theta v)\tilde{v}'(x)^2d \theta & 0 \\
-\int_0^1p''(\tilde{v}(x)+\theta v)\tilde{v}'(x)d \theta & 0
    \end{pmatrix} \, .
\end{align*}
First, consider for every $\nu$ positive and $k$ in $\N^+$ the viscous approximations (with $U_{1,\nu}(x,t)=S(x)+U_0(x)$ for every $\nu$ positive, $t$ nonnegative and every real $x$) given by
\begin{align*}
    \partial_t(U_{k+1,\nu}-S)(x,t)=&A(U_{k,\nu}(x,t),x)\partial_x(U_{k+1,\nu}-S)(x,t) \\ &+\partial_x(B(U_{k,\nu}(x,t))\partial_x(U_{k+1,\nu}-S)(x,t)) + \nu \partial_x^2(U_{k+1,\nu}-S)(x,t) \\ &+C(U_{k,\nu}(x,t),x)(U_{k+1,\nu}(x,t)-S(x)) \, .
\end{align*}
Hence, for every $\nu$ positive and $k$ natural number, $U_{k,\nu}$ is well-defined on $\R \times [0,T_{k,\nu})$ for some $T_{k,\nu}$ positive by classical parabolic theory (it may happen that $T_{k,\nu}$ is not infinity as solutions are valued in $\R^+ \times \R$) , and for every $n$ in $\N$, it is in $C^0([0,T_{k,\nu}],H^n(\R))$. \\
Furthermore, doing standard energy estimates for the functionals 
\begin{align*}
    v \mapsto \langle \eta''(U_{k,\nu})\partial_x^{\ell}v,\partial_x^{\ell}v \rangle \, ,
\end{align*}
(where $\ell$ is in $\lbrace 0,1,2,3 \rbrace$) on $[0,T_{k,\nu})$, one can show that for any $ \inf_{k,\nu}T_{k,\nu}^*$ is positive and that for any $T$ in $(0,\inf_{k,\nu}T_{k,\nu}^*)$, $(U_{k,\nu}-S)_{k,\nu}$ is uniformly bounded in $C^0([0,T],H^3(\R))$. Then, by taking $T$ smaller, it can be proven that for any $\nu$ positive, the sequence $(U_{k,\nu}-S)_{k \in \N}$ is contractive in $C^0([0,T],L^2(\R))$ uniformly with respect to $\nu$. Finally, by compactness results, it allows to consider a family of $(U_k)_{k \in \N}$ such that $U_1(\cdot,t)=U_0$ for any $t$ and such that for any $k$ in $\N$ 
\begin{align*}
    \partial_t(U_{k+1}-S)(x,t)=&A(U_k(x,t),x)\partial_x(U_{k+1}-S)(x,t) \\ &+\partial_x(B(U_k(x,t))\partial_x(U_{k+1}-S)(x,t))\\ &+C(U_k(x,t),x)(U_{k+1}(x,t)-S(x)) \, .
\end{align*}
Finally, the contractive property of the sequence $(U_k)_k$ passes to the limit. Thus, we can consider the limit $U$ in in $C^0([0,T],L^2(\R))$. It satisfies in the distributional sense $\partial_t(U-S)=A(U)\partial_x(U-S)+C(U,\cdot)(U-S)+\partial_x(B(U)\partial_x(U-S))$ is bounded and uniformly continuous on $\R \times [0,T]$ and satisfies that $U(\cdot,x)=U_0$. Hence $U$ is in $C^0([0,T],H^2(\R))$ and we have the desired result using Sobolev's embedding theorem. \\
Finally, the uniqueness in $C^0([0,T],H^2(\R)) \cap L^{\infty}([0,T],H^3(\R))$ follows from a similar energy estimate in $L^2(\R)$.
\end{proof}
We are now ready to conclude the argument. As before, we may bound the difference:
\[
|\dot{X(t)}Y(t)+B(t)-G(t)-\dot{X(0)}Y(0)+B(0)-G(0)| \leq \epsilon,
\]
for any $\epsilon > 0$, for $0 \leq t \leq T^*$ with $T^*$ sufficiently small, by invoking Lemma \ref{approximationlemma}. Here, we use the functions (cf. \eqref{systemrepresentationkv}):
\begin{enumerate}
    \item $q_0 \left(\begin{pmatrix} v \\ h \end{pmatrix}\right)=\begin{pmatrix} p(v) \\ h \end{pmatrix}$,
    \item $\tilde{q}_0\left(\begin{pmatrix} v \\ h \end{pmatrix} \right)=v$,
    \item $q_1(u,v)=(v-u)$,
    \item $q_2(u,v)=\frac{|u_2-v_2|^2}{2}+Q(u_1|v_1)$,
    \item $q_3(u,v)=\sigma p(u_1|v_1)$,
    \item $q_4(u,v)=\frac{1}{2}(p(u_1)-p(v_1))^2$,
    \item $q_5(u,v)=p(u_1)-p(v_1)$,
    \item $q_6(u,v)=\frac{1}{2\sigma}(p(u_1)-p(v_1))^2-\frac{\sigma}{2}\left(u_2-v_2-\frac{p(u_1)-p(v_1)}{\sigma}\right)^2-\sigma Q(u_1|v_1)$.
\end{enumerate}
So, taking the limit as in the scalar case, we obtain the result. 
\appendix
\section{Proof of Lemma \ref{approximationlemma}}
Here, we give a proof for the general lemma that is used to show that the time derivative of the relative entropy is continuous uniformly in $n$ for a sequence of $C^1$ shifts approximating a Lipschitz one. This is a key lemma needed to show our main results for Lipschitz shifts rather than $C^1$ shifts. 
\begin{proof}
It suffices to show the supremum over $X,t$ may be made arbitrarily small in the difference of each pair of integrals involving $q_i$ independently. For the terms involving $q_1, q_2$, we check: 
\begin{align*}
&\Bigl|\dot{X}(t)\left(\int_{\R}a\left(q_0(\tilde{w})\right)_x \cdot q_1(w^{X(t)}(t),\tilde{w})\diff x+\int_{\R}a'q_2(w^{X(t)}(t),\tilde{w})\diff x\right) \\
&-\dot{X}(0)\left(\int_{\R}a\left(q_0(\tilde{w})\right)_x \cdot q_1(w(0),\tilde{w})\diff x-\int_{\R}a'q_2(w(0),\tilde{w})\diff x\right)\Bigr|\\
=& \Bigl|\dot{X}(t)\Bigl(\int_{\R}a\left(q_0(\tilde{w})\right)_xq_1(w^{X(t)}(t),\tilde{w})\diff x-\int_{\R}a\left(q_0(\tilde{w})\right)_xq_1(w(0),\tilde{w})\diff x \\
&+\int_{\R}a'q_2(w^{X(t)}(t),\tilde{w})\diff x-\int_{\R}a'q_2(w(0),\tilde{w})\diff x\Bigr)\Bigr|  \\
\leq& \|\dot{X}\|_\infty \Bigl(\|\tilde{a}\|_\infty\|\left(q_0(\tilde{w})\right)_x\|_2\|\nabla q_1\|_{\infty}\|w^{X(t)}(t)-w(0)\|_2 \\
&+\|\tilde{a}'\|_\infty\|\tilde{w}'\|_2\|\nabla q_2\|_\infty\|w^{X(t)}(t)-w(0)\|_2\Bigr).
\end{align*}
Further, using Lemma $\ref{shifttranslate}$, we obtain:
\[
\|w^{X(t)}(t)-w(0)\|_2 \leq \|w_x(t)\|_2|X(t)|+\|w(t)-w(0)\|_2.
\]
For $0 \leq t \leq T$, $\|w_x(t)\|_2$ is bounded, and $|X(t)| \leq Kt$. Further, $\|w(t)-w(0)\|_2$ goes to $0$ as $t$ goes to $0$ due to the strong continuity into $L^2$. So, the difference may be made arbitrarily small (uniformly in $X$) for $t$ sufficiently small. 
\par For the terms involving $q_3$, we check:
\begin{align*}
&\bigl|\int_{\R}a\left(\tilde{q}_0(\tilde{w})\right)_xq_3(w^{X(t)}(t),\tilde{w})\diff x-\int_{\R}a\left(\tilde{q}_0(\tilde{w})\right)_xq_3(w(0),\tilde{w})\diff x\bigr| \\
&\leq \|\tilde{a}\|_\infty \|\left(\tilde{q}_0(\tilde{w})\right)_x\|_2\|\nabla q_3\|_\infty\|w^{X(t)}(t)-w(0)\|_2.
\end{align*}
This may be made small uniformly in $X$ in the same way as above.
\par For the terms involving $q_4$, we check:
\begin{align*}
&\bigl|\int_{\R}a''q_4(w^{X(t)}(t),\tilde{w})\diff x-\int_{\R}a''q_4(w(0),\tilde{w})\diff x\bigr| \\
&\leq \|a''\|_2\|\nabla q_4\|_\infty \|w^{X(t)}(t)-w(0)\|_2.
\end{align*}
This may be made small uniformly in $X$ as long as $a'' \in L^2$. We check:
\[
a''=\tilde{a}''(\tilde{w}_1)(\tilde{w}_1')^2+\tilde{a}'(\tilde{w}_1)\tilde{w}_1''.
\]
So, $a'' \in L^2$ because $\tilde{w}_1'' \in L^2$, and $\tilde{w}_1' \in L^4$.
\par For the terms involving $q_5$, we check:
\begin{align*}
&\bigl|\int_{\R}a\left(q_5(w^{X(t)}(t),\tilde{w})\right)_x^2\diff x-\int_{\R}a\bigl(q_5(w(0),\tilde{w})\bigr)_x^2\diff x\bigr| \\
&=\bigl|\int_{\R}a\left((q^5(w_1^{X(t)}(t))-q^5(\tilde{w}_1))_x^2-(q^5(w_1(0))-q^5(\tilde{w}_1))_x^2\right)\diff x\bigr| \\
&=\bigl|\int_{\R}a\left(A^2+2AB\right)\diff x\bigr|,
\end{align*}
where:
\begin{align*}
&A(t)=(q^{5})'(w_1^{X(t)}(t))(w_{1})_x^{X(t)}(t)-(q^{5})'(w_1(0))(w_{1})_x(0), \\
&B=(q^{5})'(w_1(0))(w_{1})_x(0)-(q^{5})'(\tilde{w}_1)\tilde{w}_1'.
\end{align*}
So, we obtain:
\[
\bigl|\int_{\R}a\left(A^2+2AB\right)\diff x\bigr| \leq \|\tilde{a}\|_\infty (\|A\|_2^2+2\|A\|_2\|B\|_2).
\]
We will show that $\|A\|_2$ goes to $0$ as $t$ goes to $0$ uniformly in $X$, and $\|B\|_2$ is bounded. For $A$:
\begin{align*}
&\int_{\R}|(q^{5})'(w_1^{X(t)}(t))(w_1)_x^{X(t)}(t)-(q^{5})'(w_1(0))(w_1)_x(0)|^2 \diff x \\
\leq& \int_{\R}\left(|(q^{5})'(w_1^{X(t)}(t))||(w_1)_x^{X(t)}(t)-(w_1)_x(0)| \right. \\
&\left. +\|(q^{5})''\|_\infty |w_1^{X(t)}(t)-w_1(0)\|(w_1)_x(0)|\right)^2\diff x \\
 \leq& 2\|(q^{5})'\|_\infty^2 \|(w_1)_x^{X(t)}(t)-(w_1)_x(0)\|_2^2 \\
& +2\|(q^{5})''\|_\infty^2\|(w_1)_x(0)\|_\infty^2\|w_1^{X(t)}(t)-w_1(0)\|_2^2.
\end{align*}
Both $\|(w_1)_x^{X(t)}(t)-(w_1)_x(0)\|_2$ and $\|w_1^{X(t)}(t)-w_1(0)\|_2$ converge to 0 uniformly in $X$ as $t$ goes $0$ by Lemma \ref{shifttranslate}. For $B$:
\begin{align*}
&\int_{\R}|(q^{5})'(w_1(0))(w_1)_x(0)-(q^{5})'(\tilde{w}_1)\tilde{w}_1'|^2dx \\
&\leq 2\|(q^5)'\|_\infty^2\|(w_1)_x(0)-\tilde{w}_1'\|_2^2+2\|(q^5)''\|_\infty^2\|\tilde{w}_1'\|_\infty^2\|w_1(0)-\tilde{w}_1\|_2^2.
\end{align*}
Finally, the terms involving $q_6$ are handled identically to the terms involving $q_4$, as $a' \in L^2$.
\end{proof}

\bibliographystyle{plain}
\bibliography{refs}
\end{document}